\documentclass[a4paper,10pt]{article}
\usepackage{amssymb}
\usepackage{latexsym}
\usepackage{amsmath}
\usepackage{amsthm}
\usepackage{amsfonts}
\usepackage{amssymb}
\usepackage{amsmath,amscd}
\usepackage{graphicx}

\newtheorem{Th}{Theorem}
\newtheorem{Prop}{Proposition}
\newtheorem{Co}{Corollary}
\newtheorem{Lm}{Lemma}

\newtheorem{Dfi}{Definition}
\newtheorem{Rm}{Remark}

\newcommand{\be}{\begin{equation}}
\newcommand{\ee}{\end{equation}}
\newcommand{\bes}{\begin{equation*}}
\newcommand{\ees}{\end{equation*}}
\newcommand{\vertiii}[1]{{\left\vert\kern-0.25ex\left\vert\kern-0.25ex\left\vert #1 
    \right\vert\kern-0.25ex\right\vert\kern-0.25ex\right\vert}}
\newcommand{\R}{\mathbb{R}}

\newcommand{\reset}{\setcounter{equation}{0}\setcounter{Th}{0}\setcounter{Prop}{0}\setcounter{Co}{0}
\setcounter{Lm}{0}\setcounter{Rm}{0}}

\def\lf{\left}
\def\rg{\right}

\def\al{\alpha}

\def\ep{\varepsilon}

\def\ds{\displaystyle}

\def\Om{\Omega}
\def\om{\omega}

\begin{document}

\title{3-Commutators Revisited}

\author{ Francesca Da Lio and Tristan Rivi\`ere\footnote{Department of Mathematics, ETH Zentrum,
CH-8093 Z\"urich, Switzerland.}}

\maketitle

{\bf Abstract :}{\it We present a class of Pseudo-differential elliptic systems with anti-self-dual potentials on ${\R}$ satisfying compensation phenomena similar to the ones discovered in \cite{Riv1} for elliptic systems with anti-symmetric potentials. These compensation phenomena are based on new ``multi-commutator'' structures generalizing the 3-commtators introduced by the authors in \cite{DLR1}. }

\par\medskip
{\noindent{\small { \bf  MSC 2000.} 35R09, 45K05, 42B37, 47B34, 47B47, 47B38}

 \tableofcontents 
\section{Introduction}
In the paper \cite{DLR2} the authors discovered the following  compensation phenomenon: 
\begin{equation}\left\{\begin{array}{c}
 \mbox{ if $\Omega\in L^2_{loc}(\R,so(m)),$ $v\in L^2_{loc}(\R,\R^m)$
and $f\in L^p_{loc}(\R,\R^m)$ ($2>p\ge1$) satisfy}\\[5mm]
(-\Delta)^{1/4} v=\Omega v+f,\label{nonlocalSch}\\[5mm]
 \mbox{then ~~$(-\Delta)^{1/4} v\in L^p_{loc}(\R)$}.\end{array}\right.\end{equation}

This result is central in the regularity theory of $\frac{1}{2}$-harmonic map. It came after a similar theorem for local   elliptic Schr\"odinger type systems with an antisymmetric potential \cite{Riv1}.\par
These phenomena are based on the existence of special linear operator satisfying ``better" integrability properties due to compensation.\par
Such an operator is for instance given by the so-called  $3$-commutator:\par
\begin{equation}\label{3term}
{\cal{T}}\colon(v,Q)\mapsto (-\Delta)^{1/4} (Qv)-Q\, (-\Delta)^{1/4} v +(-\Delta)^{1/4} (Q) v.
\end{equation}
  It  is proved in \cite{DLR1} that ${\mathcal T}\colon L^2(\R)\times \dot{H}^{1/2}(\R)\to H^{-1/2}(\R)$ and 
\begin{equation}\label{firstope}
\|{\cal{T}}(v,Q)\|_{H^{-1/2}(\R)}\lesssim \|Q\|_{\dot{H}^{1/2}(\R)}\|v\|_{L^{2}(\R)}.
\end{equation}
 
\par
The operator ${\cal{T}}$  appears as a natural replacement for $1/2$-harmonic maps of the existing  Jacobian structures for harmonic maps
into manifolds (see developments on that topic in \cite{DLS2}.)\par
The   estimate \eqref{firstope} has been originally proved using Littlewood-Paley  dyadic  decomposition (see an alternative proof in \cite{LS}).\par
In the  present work we are going  to generalize the $3$-commutators \eqref{3term} to ``multi-commutators" enjoying similar compensation phenomena.
These commutators will be useful to deduce  regularity results for
integro-differential elliptic systems of the form
\begin{equation}\label{systemnew-0}
(-\Delta)^{1/4} v=\int_{\R} H(x,y) v(y) dy +f(x),
\end{equation}
where  $H$  satisfies suitable conditions that we are now specifying. \par
 \par

We introduce the following  {\bf Besov type spaces of Schwarz Kernels} : for $s\in\R$, $1<p<+\infty$, $1\le q< +\infty$,   
 we denote by    $A_{p,q}^s(\R^{2n},M_m({\R}))$ the following space
  \be
  \label{A1-0}
 \ds A_{p,q}^s(\R^{2n},M_m({\R}))=\lf\{K\in L^{1}_{loc}(\R^{2n},M_m({\R})): ~\left(\int_{\R^n}|h|^{n-qs}\|K(\cdot,\cdot+h) )\|^q_{L^{p}(\R^n)}\,dh\right)^{1/q}<+\infty\rg\}
 \ee
 Our main result is the following
 \begin{Th}
 \label{th-intro.1}
 Let $K\in L^1_{loc}(\R^{2},M_m({\R}))$ and $0<\sigma<1/2$ such that $(-\Delta)^{\sigma/2}K\in A^{-\sigma}_{2,2}({\R}^{2},M_m({\R}))$ 
 \footnote{We recall that for $f\in {\mathcal{S}}(\R),$ ( ${\mathcal{S}}(\R)$ is the space of Schwarz functions) $(-\Delta)^{\sigma/2}f(x):= pv\int_{\R}\frac{f(x)-f(y)}{|x-y|^{1+\sigma}} dy$} that is to say
 \be
 \label{assker}
 \int_{\R}\int_{\R} |h|^{1+2\sigma}\ |(-\Delta_x)^{\sigma/2}K(x+h,x)|^2\ dx\,dh<+\infty\quad,
 \ee
 and let $\om\in L^1_{loc}({\R},M_m({\R}))$ such that
 \be
 \label{assomega}
 \om(x)-\int_{\R}K^t(x,y)\, dy\ \in \ L^2({\R},so(m))
 \ee
 and
 \be
 \label{assKe}
 K(x,y)=-\,K^t(y,x)\quad \mbox{for  a.e.}\,  (x,y)\quad.
 \ee
 Then for any $v\in L^2_{loc}({\R},{\R}^m)$ and any $f\in L^p_{loc}({\R},{\R}^m)$ where $1\le p<2$ satisfying
 \be
 \label{systemnew-1}
 (-\Delta)^{1/4} v=\int_{\R} H(x,y)\ v(y) dy +f(x)\quad,
 \ee
where $H(x,y):=K(x,y) +\om(x)\, \delta_{x=y}$ we have
 \be
 \label{intro-1}
 (-\Delta)^{1/4} v\in L^p_{loc}({\R})\quad.
 \ee
 \hfill $\Box$
 \end{Th}
\begin{Rm}
\label{rm-intro.1}
 Operators whose kernel satisfy the  condition  \eqref{assKe} are  {\bf anti-self-adjoint } operators. We believe that the above result
 should be generalised to more general { anti-self adjoint } operators with the ad-hoc mapping properties corresponding to the membership of $(-\Delta)^{\sigma/2}K$ in $A^{-\sigma}_{2,2}({\R}^{2},M_m({\R}))$ for the operator
 \[
 v\ \longrightarrow \int_{\R}K(x,y)\ v(y)dy\quad.
 \]
 Hence the {\bf anti-symmetry} condition which was the key notion in the original work \cite{Riv1} should be somehow substituted by the more general  notion 
 of {\bf anti-self-adjointness}.\hfill $\Box$
 \end{Rm}
 \begin{Rm}
 \label{rm-intro.1a}
 It would be interesting to study the possibility to extend or not the previous theorem to the limiting case $\sigma=0$ and the assumption
 \[
 \int_{\R}\int_{\R} |x-y|\ |K(x,y)|^2\ dx\,dy<+\infty
 \]
 instead of (\ref{assker}).
 \hfill $\Box$
 \end{Rm}
We need first to clarify  the  meaning to the integral $\int_{\R}H(x,y)v(y)dy$ under the assumptions  (\ref{assker}), \eqref{assomega} and \eqref{assK} and $v$ in $L^2$ .
This justification  is based on the notion of {\bf abstract multicommutator} that we are introducing now.
It is not a-priori clear that under the above assumptions on $K$ one has that
 $H(x,y)\in L^2(\R)$ for a.e. $x\in\R$. Nevertheless we can give $\int_{\R}H(x,y)v(y)dy$    a meaning in 
 a distributional sense.\par
  Precisely we introduce the following definition.
  \begin{Dfi}
  \label{df-intro.1} {\bf[Abstract multi-commutators]}
  Let $K\in L^1_{loc}({\R}^{2n},M_m({\R}))$ satisfying the anti-self-dual condition
  \be
  \label{anti-self-dual}
  K(x,y)=-\,K^t(y,x)\quad \mbox{ for a. e. } (x,y)\in {\R}^{2n}\quad,
  \ee
  then the operator given by
  \be
\label{I.3}
{\mathcal T}_K(v)(x):=\int_{\R}\lf[K^t(x,y)\, v(x)+K(x,y)\ v(y)\rg]\ dy
\ee
is called an  { abstract multi-commutator}.\hfill $\Box$
\end{Dfi}
Such an operator enjoys the following integrability by compensation property\footnote{Similar property hold in higher dimension obviously} in ${\R}$ 
 \begin{Lm}
\label{lm-I.1} {\bf [Compensation for multi-commutators]} Under the above notations let $K$ be a Schwartz Kernel satisfying the anti-self-dual condition (\ref{anti-self-dual}) one has for any $r>1$, $p>r'$, $q\ge 2$ and $\sigma>0$
\be
\label{I.4}
\|{\mathcal T}_K(v)\|_{\dot B^{-(2/q-1+\sigma)}_{rp/(p+r),q'}({\R})}\le C\ \|K\|_{A^{-\sigma}_{p,q}({\R}^2)}\ \|v\|_{L^r({\R})}\quad.
\ee
where $\dot B^{s}_{p,q}({\R})$ denotes the usual homogeneous Besov spaces~\footnote{For  $s\in \R$, $1<p<+\infty$ and $1\le q < +\infty$ we also denote by $\dot B^s_{p,q}(\R^n)$  the homogeneous   Besov spaces  given by:
\begin{eqnarray}
\ds\dot B^s_{p,q}(\R^n)&=&\lf\{f\in L^{1}_{loc}(\R^n): ~\left(\int_{\R^n}|h|^{-n-sq}\|f(x+h)-f(x)\|^q_{L^{p}(\R^n)}\,dh\right)^{1/q}<+\infty\rg\}\label{Besov1}\\
 \end{eqnarray}
}.\hfill $\Box$
\end{Lm}
In order now to justify the integral in the r.h.s. of (\ref{systemnew}) we can then write
 \begin{eqnarray}
 \int_{\R} H(x,y) v(y)dy&=&{\mathcal T}_K(v)(x)+\left(\omega(x)-\int_{\R}K^t(x,y) dy\right)v(x).
 \end{eqnarray}
 Besides proving the main theorem~\ref{th-intro.1}, the goal of the paper is to illustrate the relative easiness  to produce multi-commutators.
 
 \medskip

 For instance the $3$-term commutator \eqref{3term} is  a particular example of operator ${\mathcal T}_K$ where  $K:=K_{d^{1/2}Q}$ is the Schwartz Kernel associated
 to the operator
 \[
 d^{1/2}Q:=Q\circ(-\Delta)^{1/4}-(-\Delta)^{1/4}\circ Q\quad .
 \]
It is given explicitly by
\begin{equation}\label{KT}
K_{d^{1/2}Q}(x,y):= \frac{Q(x)-Q(y)}{|x-y|^{3/2}}\quad.
\end{equation}
One has obviously
\be
 \label{II.2a}
 K_{d^{1/2}Q}(x,y)=-\,K^t_{d^{1/2}Q}(y,x)
 \ee
 and a direct computation gives
\be
\label{II.3}
\begin{array}{l}
\ds {\mathcal T}_{ K_{d^{1/2}Q}}(v)=\lf[Q\circ(-\Delta)^{1/4}-(-\Delta)^{1/4}\circ Q- (-\Delta)^{1/4}Q\rg]v=\int_{\R}\frac{Q(y)-Q(x)}{|x-y|^{3/2}} \ [v(y)+v(x)]dy \\[5mm]
\ds\quad=\int_{R} \lf[K_{d^{1/2}Q}(x,y)\, v(y)-K_{d^{1/2}Q}(y,x) v(x)\rg]dy
\end{array}
\ee
Observe that for any $0\le \sigma\le 1/2$ one has
\be
\label{intro-a}
\|(-\Delta)^{\sigma/2}K_{d^{1/2}Q}\|^2_{A^{-\sigma}_{2,2}}= \int_{\R}\int_{\R} |h|^{1+2\sigma}\ |(-\Delta_x)^{\sigma/2}K_{d^{1/2}Q}(x+h,x)|^2\ dx\,dh\simeq \|(-\Delta)^{\sigma/2}Q\|^2_{\dot{B}^{1/2-\sigma}_{2,2}}\simeq \|Q\|^2_{\dot{B}^{1/2}_{2,2}}
\ee
In particular Lemma \ref{lm-I.1} implies the compensation phenomena observed in \cite{DLR1} for  \eqref{3term}. Indeed, Sobolev embedding implies the continuity of the map
\be
\label{embker}
(-\Delta)^{-\sigma/2} \ :\  A^{-\sigma}_{2,2}({\R})\ \longrightarrow\ A^{-\sigma}_{2/(1-2\sigma), 2}({\R})
\ee
We can then apply  lemma \ref{lm-I.1} $p=4$, $q=2$ and $\sigma=1/4$ and using the facts that 
\[
\dot{B}^{-1/4}_{4/3,2}({\R})=(\dot{B}^{1/4}_{4,2}({\R})'\hookrightarrow H^{-1/2}({\R})=\dot{B}^{-1/2}_{2,2}({\R})=(\dot{B}^{1/2}_{2,2}{\R})'.
\]
(see e.g. \cite{RS}) we recover (\ref{firstope}).
\begin{Rm}
\label{rm-intro.2}
In \cite{MS} the authors have discovered a compensation phenomenon for general $L^2$ Kernels which are $div_{1/2}-$free and contracted with $K_{d^{1/2}Q}$ 
for any $Q\in \dot{H}^{1/2}({\R})$. While this compensation phenomenon, which is a fractional version of Wente-Coifman-Lions-Meyer-Semmes compensation, is obviously of similar nature to the one given in lemma~\ref{lm-I.1}, they seem however not to be completely ``isomorphic'' to each other.\hfill $\Box$
\end{Rm}
Another elementary but useful observation is the {stability} of the  property \eqref{assK}     with respect to the 
 adjoint multiplication by $P\in L^{\infty}(\R,M_{m}({\R}))$ .
    Precisely we have  
\begin{Lm}\label{stabKernel}   {\bf[Stability of multi-commutators by adjoint multiplication.]} Let  $K$ satisfy  \eqref{assK} and $P\in L^{\infty}(\R,M_{m}({\R}))$. Then  the
new kernel $$G(x,y):=P(x)\,K(x,y)\,P^t(y)$$  satisfies  \eqref{assK}   and defines a new multi-commutator. In particular 
 $${\mathcal T}_{G}(v)(x):=\int_{\R}\lf[G^t(x,y)\, v(x)+G(x,y)\ v(y)\rg]\ dy$$
 satisfies the compensation lemma \ref{lm-I.1}. \hfill $\Box$
 \end{Lm}
 \par
As far as the stability  with respect to the composition   with a pseudo-differential  operator of order zero is concerned in the present paper we focus our attention to the composition 
between the Riesz operator ${\mathfrak R}$ and $d^{1/2}Q$ for any $Q\in \dot{H}^{1/2}({\R},\mbox{Sym}_{m}({\R}))$  where $\mbox{Sym}_{m}({\R})$ denotes the space of square $m$ by $m$ real matrices.
\begin{Lm}
\label{lm-II.1}{\bf[Genarating a multicommutator from ${\mathfrak R}\circ d^{1/2}Q$]}
Let $Q\in \dot{H}^{1/2}({\R},\mbox{Sym}_{m}({\R}))$ then the Schwartz Kernel ${R}^Q(x,y)$ of the following operator
\[
{\mathcal R}^Q:={\mathfrak R}\circ  d^{1/2}Q-   d^{1/2}Q\circ{\mathfrak R}-  {\mathfrak R}\circ(-\Delta)^{1/4}Q-(-\Delta)^{1/4}Q\circ{\mathfrak R}
\] satisfies for any $p\ge 2$ 
\be
\label{II.14}
\| {R}^Q(x,y)\|_{A_{p,2}^{-1/2+1/p}({\R}^2)} \le C_p\ \|Q\|_{\dot{H}^{1/2}({\R})}\quad.
\ee
\hfill $\Box$
\end{Lm}
Hence we deduce the following Corollary
\begin{Co}
\label{co-II.1intr}
Let $Q\in \dot{H}^{1/2}({\R},\mbox{Sym}_m({\R}))$ then the operator given by
\[
{\mathcal T}_{R^Q}:={\mathfrak R}\circ  d^{1/2}Q-   d^{1/2}Q\circ{\mathfrak R}-  {\mathfrak R}\circ(-\Delta)^{1/4}Q-(-\Delta)^{1/4}Q\circ{\mathfrak R}
-2\,{\mathfrak R}\circ(-\Delta)^{1/4}Q 
\] 
is a multi-commutator to which the compensation lemma~\ref{lm-I.1} applies.
\hfill $\Box$
\end{Co}
We can move on in complexity and consider the composition of $d^{1/2}Q$ on the right and on the left by $\mathfrak R$. Precisely we have
\begin{Lm}
\label{lm-d^{1/2}Rintro}{\bf[Genarating a multicommutator from ${\mathfrak R}\circ d^{1/2}Q\circ{\mathfrak R}$.]}
Let $Q\in \dot{H}^{1/2}({\R},\mbox{Sym}_m({\R}))$ then the Schwartz Kernel ${S}^Q(x,y)$ of the following operator
\[
{\mathcal S}^Q:={\mathfrak R}\circ d^{1/2}Q\circ {\mathfrak R}+ {\mathfrak R}[(-\Delta)^{1/4}Q]\circ {\mathfrak R}+{\mathfrak R}\circ {\mathfrak R}[(-\Delta)^{1/4}Q]
\]
satisfies for any $p>2$
\be
\label{II-f-30bis}
\lf\|S^Q\rg\|_{A^{-1/2+1/p}_{p,2}}\le  C_p\ \|Q\|_{\dot{H}^{1/2}({\R})}
\ee\quad.
\hfill $\Box$
\end{Lm}
Hence we deduce the following Corollary
\begin{Co}
\label{co-II.2intr}
Let $Q\in \dot{H}^{1/2}({\R},\mbox{Sym}_m({\R}))$
 then the operator given by
\be
\label{intro-10}
{\cal{T}}_{{\cal{S}}^Q}:={\mathfrak R}\circ d^{1/2}Q\circ {\mathfrak R}+ {\mathfrak R}[(-\Delta)^{1/4}Q]\circ {\mathfrak R}+{\mathfrak R}\circ {\mathfrak R}[(-\Delta)^{1/4}Q]- \,(-\Delta)^{1/4}Q\quad.
\ee
is a multi-commutator to which the compensation lemma~\ref{lm-I.1} applies.
\hfill $\Box$
\end{Co}

It would be interesting to explore the general rule for producing multi-commutators starting from the composition of a general anti-self-dual kernel $K$  with the Riesz potential in the same spirit as Lemma \ref{lm-II.1} and Lemma \ref{lm-d^{1/2}Rintro} where the ``starting'' kernel is $K_{d^{1/2}Q}$. 
We do believe that the somehow lengthy computations from part IV, whose expositions in this work are mostly meant to be illustrative, may contain some ``genericity'' and should be inspiring for such a later purpose.

Finally, it would be interesting to establish some  relation between the membership for a Schwartz Kernel $K$  in a Besov space of Kernels $A^\sigma_{p,q}$ and  mapping properties of ${\mathcal T}_K$.
 
 \noindent{\bf Acknowledgments :} {\it A large part of the present work has been conceived while the two authors were visiting the Institute for Advanced Studies
in Princeton. They are very grateful to the IAS for the hospitality.}

 \section{Proofs of Lemma \ref{lm-I.1} and some other properties}

In this section we will prove Lemma \ref{lm-I.1}  and some other  properties of the operator ${\mathcal T}_K$ defined in  \eqref{I.3}.

\noindent{\bf Proof of lemma~\ref{lm-I.1}.} We prove the Lemma in the case $r=2$. We will prove it by duality. Let $\varphi\in \dot{B}^{s}_{2p/(p-2),q'}({\R},{\R}^m)$,
where $s=\sigma+2/q-1.$ We have
\be
\label{I.5}
\|\varphi\|_ {\dot{B}^{s}_{2p/(p-2),q'}({\R},{\R}^m)}\simeq\left(\int_{\R}\frac{dh}{h^{1+q' s}}\,\|\varphi(\cdot+h)-\varphi(\cdot)\|^{q'}_{L^{2p/(p-2)}({\R})}\right)^{1/q'}\quad.
\ee
We compute
\be
\label{I.6}
\begin{array}{l}
\ds\lf<\varphi,{\mathcal T}_K(v)\rg>_{\dot{B}^{s}_{2p/(p-2),q'},\dot{B}^{-s}_{2p/(p+2),q}}:=\int_{\R}\int_{\R}\varphi(x)\left(K^t(x,y))\, v(x)+K(x,y)\ v(y)\rg)\ dx\ dy\\[5mm]
\ds\quad\quad\quad\quad\quad=\int_{\R}\int_{\R}\varphi(y)\lf(K^t(y,x)\, v(y)+K(y,x)\ v(x)\rg)\ dx\ dy\\[5mm]
\ds\quad\quad\quad\quad\quad=-\int_{\R}\int_{\R}\varphi(y)(K(x,y)\, v(y)\ dx\ dy+\int_{\R}\int_{\R}\varphi(x)K(x,y)\ v(y)\ dx\ dy\\[5mm]
\ds\quad\quad\quad\quad\quad=\int_{\R}\int_{\R}\lf(\varphi(x)-\varphi(y)\rg)K(x,y)v(y)\ dx\ dy
\end{array}
\ee
Hence we deduce
\be
\label{I.7}
\lf<\varphi,{\mathcal T}_K(v)\rg>_{\dot{B}^{s}_{2p/(p-2),q'},\dot{B}^{-s}_{2p/(p+2),q}}\le \ \int_{\R}\int_{\R}\ |K(x,y)|\ |\varphi(x)-\varphi(y)|\ |v(y)|\ dx dy
\ee
H\"older inequality gives
\be
\label{I.8}
\lf<\varphi,{\mathcal T}_K(v)\rg>_{\dot{B}^{s}_{2p/(p-2),q'},\dot{B}^{-s}_{2p/(p+2),q}}\le 2\, \|v\|_{L^2}\ \ \int_{\R}\, \|K(\cdot,\cdot+h)\|_{L^p(\R)}\ \ \|\varphi(\cdot+h)-\varphi(\cdot)\|_{2p/(p-2)}\ dh
\ee
Using Cauchy Schwartz this time one has
\begin{eqnarray}
\label{I.9}
&&\lf<\varphi,{\mathcal T}_K(v)\rg>_{\dot{B}^{s}_{2p/(p-2),q'},\dot{B}^{-s}_{2p/(p+2),q}}\nonumber\\[5mm]&&~~~\le 2\, \|v\|_{L^2} \left( \int_{\R}|h|^{1+\sigma q} \|K(\cdot,\cdot+h)\|^q_{L^p(\R)}\,dh\right)^{1/q}
\left( \int_{\R}|h|^{-q'/q-\sigma q'}\|\varphi(\cdot+h)-\varphi(\cdot)\|^{q'}_{2p/(p-2)}\,dh\right)^{1/q'}\nonumber\\[5mm]
&&~~~=\|v\|_{L^2} \left( \int_{\R}|h|^{1+\sigma q} \|K(\cdot,\cdot+h)\|^q_{L^p(\R)}\,dh\right)^{1/q}
\left( \int_{\R}|h|^{-1-sq'}\|\varphi(\cdot+h)-\varphi(\cdot)\|^{q'}_{2p/(p-2)}\,dh\right)^{1/q'}\nonumber\\[5mm]
&&~~~=\|v\|_{L^2}\ \|K\|_{A^{-\sigma}_{p,q}}\ \|\varphi\|_{\dot{B}^{s}_{2p/(p-2),q'}}.
\end{eqnarray}

Which proves the lemma.\hfill $\Box$

Next we will show some stability property of the operator ${\mathcal T}_K$ in the case $n=1$ with respect to the multiplication by $P \in L^{\infty}(\R,M_m({\R}))$ .

\medskip

\noindent{\bf Proof of Lemma \ref{stabKernel}.}

  We have
$$
 G^t(x,y)=P(y)K^t(x,y)P^t(x)=-G(y,x).$$

$G(x,y)$ is in $   A^{\sigma}_{p,q}(\R^2)$ since

 \begin{equation}
\label{I.G}
\|G\|_{A^{\sigma}_{p,q}(\R^2)}:=\lf[\int_{\R}h^{1-q\sigma}\|P(\cdot)K(\cdot,\cdot+h)P^t(\cdot+h)\|^q_{L^p}\ dh\rg]^{1/q}\lesssim
\|P\|^2_{L^{\infty}(\R)}\|K\|_{A^{\sigma}_{p,q}(\R^2)}<+\infty.~~\Box
\end{equation}

Next we show another property of a kernel $K$ such that $(-\Delta)^{\sigma/2}K\in  A^{-\sigma}_{2,2}(\R^2)$. \par To this purpose we extend naturally the Besov space of 
Schwartz Kernels (\ref{A1}) to the {\bf Lorentz-Besov Space of Schwartz Kernels}. For  $s\in\R$, $1<p<$, $1\le q\le +\infty$, and $1\le r\le +\infty$ 
 we denote by    $A_{(p,r),q}^s(\R^{2n},M_m({\R}))$ the following space \be
 \ds A_{(p,r),q}^s(\R^{2n})=\lf\{K\in L^{1}_{loc}(\R^{2n},M_m({\R})): ~\left(\int_{\R^n}|h|^{n-qs}\|K(\cdot,\cdot+h) )\|^q_{L^{p,r}(\R^n)}\,dh\right)^{1/q}<+\infty\rg\}\label{Besov-Schwartz}\\
 \ee
where $L^{p,r}({\R}^n)$ denote the usual Lorentz spaces (see \cite{Gra1}). \par
We  now prove the following result.
\begin{Lm}\label{LemmamultP}
Let $0<\sigma<1/2$ and $K\in L^1_{loc}({\R}^2)$ such that $(-\Delta)^{\sigma/2}K\in A^{-\sigma}_{2,2}(\R^2)$ and let $P\in \dot{H}^{1/2}(\R^2)$. Then there exists a constant $C_\sigma$ depending only on $\sigma$ such that  
 \begin{equation}\label{estL21}
\lf \| \int_{\R}(P(x)-P(y))K^t(x,y)\,dy\rg\|_{L^{2,1}(\R)}\le\ C_\sigma\  \|(-\Delta)^{\sigma/2}K\|_{A^{-\sigma}_{2,2}({\R}^2)}\ \|P\|_{\dot{H}^{1/2}(\R)}\quad.
 \end{equation}
 \hfill $\Box$
 \end{Lm}
\noindent {\bf Proof of Lemma \ref{LemmamultP}.}
 We first recall the improved Sobolev embedding for any $\sigma>0$ (see \cite{ST})
 \[
   \dot{H}^{1/2}(\R)\hookrightarrow \dot B^{\sigma}_{(\sigma^{-1},2), 2}({\R})
 \]
 where
 \[
 \dot B^s_{(p,r),q}(\R^n)=\lf\{f\in L^{1}_{loc}(\R^n): ~\left(\int_{\R^n}|h|^{-n-sq}\|f(x+h)-f(x)\|^q_{L^{p,r}(\R^n)}\,dh\right)^{1/q}<+\infty\rg\}
 \]
 Then
 \begin{eqnarray*}
 \left(\int_{\R}|h|^{-1-2\sigma}\|P(\cdot+h)-P(\cdot)\|^2_{L^{(\sigma^{-1},2)}(\R)}\,dh\right)^{1/2}&\le& C_\sigma\, \left(\int_{\R }|h|^{-2}\| P(\cdot+h)- P(\cdot)\|^2_{L^{(2,2)}(\R)}\,dh\right)^{1/2}\\
 &\simeq&\|P\|_{ \dot{H}^{1/2}(\R)}.
 \end{eqnarray*}
 Observe moreover that
 \[
 C_\sigma\,  \int_{\R}\int_{\R} |h|^{1+2\sigma}\ |(-\Delta_x)^{\sigma/2}K(x+h,x)|^2\ dx\,dh\ge  \int_{\R} |h|^{1+2\sigma}\ \|K(\cdot+h,\cdot)\|^2_{L^{p,2}({\R})}
 \]
 where $p^{-1}=2^{-1}-\sigma$

 We have
 \begin{eqnarray}\label{estPK}
 &&\int_{x\in\R} v(x)\left(\int_{h\in\R}(P(x+h)-P(x))K^t(x,x+h) dh\right) dx\nonumber\\
 &&\lesssim\int_{\R}\|v\|_{L^{2,\infty}(\R)}\|(P(x)-P(x+h))K^t(x,x+h)\|_{L^{2,1}(\R)}dh\nonumber \\
 &&\lesssim \int_{\R}\|v\|_{L^{2,\infty}(\R)}h^{2^{-1}+\sigma}\|K(\cdot,\cdot+h)\|_{L^{p,2}(\R)}h^{-2^{-1}-\sigma}\|P(\cdot+h)-P(\cdot)\|_{L^{\sigma^{-1},2}(\R)}dh\nonumber\\
 &&\lesssim\|v\|_{L^{2,\infty}(\R)}\left(\int_{\R}h^{1+2\sigma}\|K(\cdot,\cdot+h)\|^2_{L^{p,2}(\R)}\right)^{1/2}\left(\int_{\R}h^{-1-2\sigma}\|P(\cdot+h)-P(\cdot)\|^2_{L^{\sigma^{-1},2}(\R)}\right)^{1/2}\nonumber\\
 &&\lesssim \|v\|_{L^{2,\infty}(\R)}\|(-\Delta)^{\sigma/2}K\|_{A^{-\sigma}_{2,2}({\R}^2)}\|P\|_{ \dot{H}^{1/2}(\R)}.
 \end{eqnarray}
 We conclude the proof of Lemma \ref{LemmamultP}.\hfill $\Box$
 \par
 \medskip
 From Lemma \ref{LemmamultP} we deduce the following result
  \begin{Prop}\label{stabpropertyintr}
 Let $0<\sigma<1/2$ and $K\in L^1_{loc}({\R}^2)$ such that $(-\Delta)^{\sigma/2}K\in A^{-\sigma}_{2,2}(\R^2)$ and let $P\in  \dot{H}^{1/2}(\R^2)$. Then
 \begin{equation}\label{stab}
 P{\mathcal T}_K(v)={\mathcal T}_G(Pv)+{\mathcal{G}}[v],
 \end{equation}
 where $$G(x,y):=P(x)K(x,y)P^t(y)\in A^{-\sigma}_{p,2}(\R^2)$$ for $p^{-1}=2^{-1}-\sigma$ and
$${\mathcal{G}}\ :\  L^{2,\infty}(\R)\ \longrightarrow \ L^{1}(\R)$$ continuously and
\[
\vertiii{{\mathcal G}}_{L^{2,\infty}\rightarrow L^{1}}   \le  C_\sigma\  \|(-\Delta)^{\sigma/2}K\|_{A^{-\sigma}_{2,2}({\R}^2)}\ \|P\|_{ \dot{H}^{1/2}(\R)}
\]
.\hfill $\Box$
  \end{Prop}
 \noindent{\bf Proof of Proposition of \ref{stabpropertyintr}.}\par
 We compute  $P{\mathcal T}_K(v)$.
\begin{eqnarray}\label{PTK}
 P(x){\mathcal T}_K(v)(x)&=&\int_{\R}\lf[P(x)K^t(x,y)\,  v(x)+P(x)K(x,y)\ v(y)\rg]\ dy\nonumber\\
 &=& \int_{\R}\left[P(y)K^t(x,y))\,P^t(x)(P(x) v(x)) +P(x)K(x,y)P^t(y)(P(y)v(y))\right]\,dy\nonumber\\&&~~+\int_{\R}[P(x)-P(y)]K(x,y))^t\, v(x)\,dy\nonumber\\
 &=&{\mathcal T}_G(Pv)(x)+\int_{\R}[P(x)-P(y)]K(x,y))^t\, v(x)\,dy.
 \end{eqnarray}
 We observe that from Lemma \ref{LemmamultP} it follows that
 $${\mathcal{G}}[v](x):=\int_{\R}[P(x)-P(y)]K^t(x,y))\, v(x)\,dy$$
 maps $ L^{2,\infty}(\R)$ into $ L^{1}(\R).$
 \hfill$\Box$
 
 \medskip
 
 Following the same argument as in the proof of lemma~\ref{lm-I.1} we establish
 \begin{Lm}
\label{lm-IIa.1} {\bf [Lorentz-Besov Compensation for multi-commutators]} Let $K$ be a Schwartz Kernel satisfying the anti-self-dual condition (\ref{anti-self-dual}). Under the above notations, for any $r>1$, $p>r'$, $q\ge 2$ and $\sigma>0$ and $t\in [1,+\infty]$
\be
\label{I.4-b}
\|{\mathcal T}_K(v)\|_{B^{-(2/q-1+\sigma)}_{(rp/(p+r),t),q'}}\le C\ \|K\|_{A^{-\sigma}_{(p,t),q}}\ \|v\|_{L^{r,\infty}}.
\ee
where $\dot B^{s}_{(p,t),q}({\R})$ denotes the usual homogeneous Lorentz-Besov spaces~\footnote{For  $s\in \R$, $1<p<+\infty$, $1\le q < +\infty$ and $t\in[1,+\infty]$ we also denote by $\dot B^s_{(p,t),q}(\R^n)$  the homogeneous   Lorentz-Besov spaces  given by:
\begin{eqnarray}
\ds\dot B^s_{(p,t),q}(\R^n)&=&\lf\{f\in L^{1}_{loc}(\R^n): ~\left(\int_{\R^n}|h|^{-n-sq}\|f(x+h)-f(x)\|^q_{L^{p,t}(\R^n)}\,dh\right)^{1/q}<+\infty\rg\}\label{Besov1b}\\
 \end{eqnarray}
}.\hfill $\Box$
\hfill $\Box$
\end{Lm}
We are now in position to prove the main theorem of the present work.
\section{Proof of theorem~\ref{th-intro.1}.}
\reset
In order to prove theorem~\ref{th-intro.1} it suffices to prove the so called ``bootstrap test'' in some space. From such a test using localization argument  from \cite{DLR1}, \cite{DLR2}, \cite{DLP} we are going to prove Morrey type decrease in the chosen space which is going to make the PDE subcritical and the regularity will follow. Precisely we are choosing the space $L^{2,\infty}$  we are going to prove the following $\epsilon-$type regularity lemma from which the main theorem~\ref{th-intro.1} can be deduced using the arguments we just described.
\begin{Th}
 \label{th-eps-reg}
 Let $K\in L^1_{loc}(\R^{2},M_m({\R}))$ and $0<\sigma<1/2$ such that $(-\Delta)^{\sigma/2}K\in A^{-\sigma}_{2,2}({\R}^{2},M_m({\R}))$  and let $\om\in L^1_{loc}({\R},M_m({\R}))$ such that
 \be
 \label{assomega-1}
 \om(x)-\int_{\R}K^t(x,y)\, dy\ \in \ L^2({\R},so(m))
 \ee
 and
 \be
 \label{assK}
 K(x,y)=-\,K^t(y,x)
 \ee
 There exists $\ep_\sigma>0$ such that, if
 \be
 \label{eps}
 \|(-\Delta)^{\sigma/2}K\|_{A^{-\sigma}_{2,2}({\R}^{2},M_m({\R}))}+\lf\|\om(x)-\int_{\R}K^t(x,y)\, dy\rg\|_{L^2({\R},so(m))}<\ep_\sigma
 \ee
 Then for any $v\in L^2({\R},{\R}^m)$  solving
 \be
 \label{systemnew}
 (-\Delta)^{1/4} v=\int_{\R} H(x,y)\ v(y) dy\quad,
 \ee
where $H(x,y):=K(x,y) +\om(x)\, \delta_{x=y}$ we have $v\equiv 0$.
 \hfill $\Box$
\end{Th}
\noindent{\bf Proof of theorem~\ref{th-eps-reg}.}
Let 
\[
\Om:=\om(x)-\int_{\R}K^t(x,y)\, dy
\]
Following \cite{DLR2} we produce $P\in \dot{H}^{1/2}({\R},SO(m))$ such that
\[
\mbox{Asym}(P^{-1}(-\Delta)^{1/4}P)=\frac{1}{2}\lf[P^{-1}(-\Delta)^{1/4}P-(-\Delta)^{1/4} (P^{-1})\, P\rg]=-\Om\quad\mbox{and}\quad\|P\|_{\dot{H}^{1/2}({\R})}\le C\ \|\Om\|_{L^2({\R})}
\]
Recall that we have
\be
\label{sym-P}
\|\mbox{Sym}(P^{-1}(-\Delta)^{1/4}P)\|_{L^{2,1}({\R})}=\lf\| \frac{1}{2}\lf[P^{-1}(-\Delta)^{1/4}P+(-\Delta)^{1/4} (P^{-1})\, P\rg]\rg\|_{L^{2,1}({\R})}\le C\ \|P\|_{\dot{H}^{1/2}({\R})}^2
\ee
Multiplying (\ref{systemnew}) by $P$, the system can be rewritten in the following form
\be
\label{syst-P}
(-\Delta)^{1/4}(Pv)={\mathcal T}(v,P)+{\mathcal T}_G(v)+{\mathcal G}[v]+P\, \mbox{Sym}(P^{-1}(-\Delta)^{1/4}P)\, v
\ee
where ${\mathcal T}$ is the 3-commutator given by (\ref{3term}) and where we used proposition~\ref{stabpropertyintr}. Recall from \cite{DLR2} that
\be
\label{3comm-2infty}
\lf\|{\mathcal T}(v,P)\rg\|_{H^{-1/2}({\R})}\le C\, \|P\|_{\dot{H}^{1/2}({\R})}\ \|v\|_{L^{2,\infty}({\R})}\quad.
\ee
Using the following precised version of (\ref{embker}) (see \cite{ST})
\[
(-\Delta)^{-\sigma/2} \ :\  A^{-\sigma}_{2,2}({\R})\ \longrightarrow\ A^{-\sigma}_{(2/(1-2\sigma),2), 2}({\R})
\]
lemma~\ref{lm-IIa.1} gives
\[
\|{\mathcal T}_G(v)\|_{B^{-\sigma}_{((1-\sigma)^{-1}, 2), 2}({\R})}\le C\,  \|(-\Delta)^{\sigma/2}K\|_{A^{-\sigma}_{2,2}({\R}^{2},M_m({\R}))}\ \|v\|_{L^{2,\infty}({\R})}
\]
Observe that for $\sigma<1/2$ we have
\[
B^{1/2}_{2,2}({\R})\ \hookrightarrow \ B^\sigma_{(\sigma,2),2}({\R})=\lf(B^{-\sigma}_{((1-\sigma)^{-1}, 2), 2}({\R})\rg)'\quad.
\]
Hence we have
\be
\label{bound-TG}
\|{\mathcal T}_G(v)\|_{\dot{H}^{-1/2}({\R})}\le C\,  \|(-\Delta)^{\sigma/2}K\|_{A^{-\sigma}_{2,2}({\R}^{2},M_m({\R}))}\ \|v\|_{L^{2,\infty}({\R})}\quad.
\ee
Combining (\ref{syst-P}), (\ref{sym-P}), (\ref{3comm-2infty}) and (\ref{bound-TG}) we obtain
\[
\|v\|_{L^{2,\infty}}\le \|(-\Delta)^{1/4}(Pv)\|_{H^{-1/2}+L^1}\le \lf[\|\Om\|_{L^2}+\|(-\Delta)^{\sigma/2}K\|_{A^{-\sigma}_{2,2}({\R}^{2},M_m({\R}))}\rg]\ \|v\|_{L^{2,\infty}({\R})}
\]
For $\ep_\sigma$ small enough in (\ref{eps}) we have $v\equiv 0$. This concludes the proof of theorem~\ref{th-eps-reg}.\hfill $\Box$

\section{Estimates of the Schwartz kernels of some commutators}
\reset
In this last section we are going to generate respectively from ${\mathcal R}\circ d^{1/2}Q$ and from ${\mathcal R}\circ d^{1/2}Q\circ {\mathcal R}$ new multi-commutators whose Schwartz Kernels satisfy the assumptions of lemma~\ref{lm-I.1}. 
Precisely we are proving lemma~\ref{lm-II.1} and lemma~\ref{lm-d^{1/2}Rintro}.

\subsection{ Producing multi-commutators from ${\mathfrak R}\circ d^{1/2}Q$.}
Let ${\mathfrak R}$ denote the Riesz transform on ${\R}$:
$$
{\mathfrak R}(v)(x):=\frac{1}{\pi}\mbox{PV} \int_{\R}\frac{v(x)-v(y)}{(x-y)} dy:=\frac{1}{\pi}\lim_{\varepsilon\to 0}\int_{|x-y|\ge \varepsilon}  \frac{v(x)-v(y)}{(x-y)} dy
$$
where PV  denotes the principal value of the integral.\par
We are going to investigate what happens if we compose the operator $d^{1/2}Q$  with the Riesz transform. We will show that we can generate from $d^{1/2}Q$ new multi-commutators that continue to satisfy Lemma~\ref{lm-I.1}. \par
To this purpose 
we compute and estimate the Schwartz Kernels of respectively
\[
R^Q:={\mathfrak R}\circ d^{1/2}Q\quad\quad\mbox{ and }\quad\quad L^Q:=d^{1/2}Q\circ {\mathfrak R}.
\]
We denote
 
\[
T_R^Q:={\mathfrak R}\circ {\mathcal{T}}_{K_{d^{1/2} Q}}\quad\quad\mbox{ and }\quad\quad T_L^Q:= {\mathcal{T}}_{K_{d^{1/2} Q}}\circ {\mathfrak R}\quad.
\]
where we recall
\be
\label{II.3-aaa}
\begin{array}{l}
\ds {\mathcal{T}}_{K_{d^{1/2} Q}}(v)=\lf[Q\circ(-\Delta)^{1/4}-(-\Delta)^{1/4}\circ Q- (-\Delta)^{1/4}Q\rg]v=\int_{\R}\frac{Q(y)-Q(x)}{|x-y|^{3/2}} \ [v(y)+v(x)]dy \\[5mm]
\ds\quad=\int_{R} \lf[K_{d^{1/2}Q}(x,y)\, v(y)-K_{d^{1/2}Q}(y,x) v(x)\rg]dy
\end{array}
\ee

We denote the corresponding Schwartz Kernels respectively by $K_{R^Q}$, $K_{L^Q}$, $K_{T_R^Q}$ and $K_{T_L^Q}$. 

Let $Q\in \dot{H}^{1/2}({\R},\mbox{Sym}_m({\R}))$. Observe that for any $v$ and $w$ in ${\mathcal S}(\R,{\R}^m)$ we have\footnote{$<\ ,\ >$ denotes the $L^2$ scalar product}, using the fact that $(-\Delta)^{1/4}$ as well as ${\mathfrak R}$ send
real functions into real functions
 \begin{eqnarray}
\label{II.4-b}
&&\ds\lf<v, (R^Q-L^Q) w\rg>= \lf<v, {\mathfrak R}\circ(Q\circ(-\Delta)^{1/4}- (-\Delta)^{1/4}Q)\,w\right>\\[5mm] &&\ds\quad\quad-\left<v,(Q\circ(-\Delta)^{1/4}-(-\Delta)^{1/4}\circ Q)\circ{\mathfrak R}\, w\rg>\nonumber\\[5mm]
&&  =\ds-\lf<{\mathfrak R}\, v,\lf(Q\circ(-\Delta)^{1/4}-(-\Delta)^{1/4}\circ Q\rg)w\rg>-\lf<\lf((-\Delta)^{1/4}\circ Q-Q\circ(-\Delta)^{1/4}\rg)\, v,{\mathfrak R}\, w\rg>\nonumber\\[5mm]
&& =\ds-\lf<\lf((-\Delta)^{1/4}\circ Q-Q\circ(-\Delta)^{1/4}\rg)\circ{\mathfrak R}\, v, w\rg>+\lf<{\mathfrak R}\circ\lf((-\Delta)^{1/4}\circ Q-Q\circ(-\Delta)^{1/4}\rg)\, v, w\rg>\nonumber\\[5mm]
&& =\ds-\lf<w, {\mathfrak R}\circ(Q\circ(-\Delta)^{1/4}-(-\Delta)^{1/4}\circ Q)-(Q\circ(-\Delta)^{1/4}-(-\Delta)^{1/4}\circ Q)\circ{\mathfrak R}\, v\rg>\nonumber\quad\\
&&=\ds-\lf<w,(R^Q-L^Q)v \rg>.
\end{eqnarray}
 The estimate \eqref{II.4-b} implies that $(R^Q-L^Q)$ is formally anti-self-dual for the $L^2$ scalar product i.e.
\be
\label{II.4-c}
(R^Q-L^Q)^\ast=-(R^Q-L^Q)\quad.
\ee
Translating this identity at the level of Schwartz Kernels give
\be
\label{II.4-d}
(K_{R^Q}-K_{L^Q})(x,y)=-(K_{R^Q}-K_{L^Q})^t(y,x)\quad.
\ee
We have
\be
\label{II.5}
(R^Qv)(x)=\frac{1}{\pi}\int_{\R}\frac{(d^{1/2}Q)v(x)-(d^{1/2}Q)v(y)}{x-y}\ dy\quad,
\ee
(the integral \eqref{II.5} is always meant in the  sense of principal value) and the calculations give
\be
\label{II.6}
K_{R^Q}(x,y)=\frac{1}{\pi}\int_{\R}\frac{dz}{x-z}\lf[\frac{Q(y)-Q(x)}{|y-x|^{3/2}}-\frac{Q(y)-Q(z)}{|z-y|^{3/2}}\rg]\quad.
\ee
On the other hand we have
\be
\label{II.KL}
K_{L^Q}(x,y)=\frac{1}{\pi}\int_{\R}\frac{dz}{y-z}\lf[\frac{Q(x)-Q(z)}{|x-z|^{3/2}} \rg]\quad.
\ee
Since $Q$ is symmetric clearly we have
\be
\label{II.6-a}
(K_{R^Q}(x,y))^t=K_{R^Q}(x,y)\quad\mbox{ and similarly }\quad(K_{L^Q}(x,y))^t=K_{L^Q}(x,y)\quad.
\ee
Hence (\ref{II.4-d}) becomes
\be
\label{II.6-b}
(K_{R^Q}-K_{L^Q})(x,y)=-(K_{R^Q}-K_{L^Q})(y,x)\quad.
\ee
 
 \medskip
 
We next prove Lemma \ref{lm-II.1}. For the simplicity of the presentation we shall restrict to the case $q=2$ and we will simply write
$A^s_p$ for $A^s_{p,2}$.  

\medskip

We set $$S_Q(x,y):=K_{R^Q}(x,y)-K_{L^Q}(x,y).$$ A priori $S_Q(x,y)$ doe not belong to the space $A_p^{-1/2+1/p}(\R^{2})$.
One has to add to the operator $R^Q-L^Q$  a suitable quantity  in oder to have a new kernel satisfying the desired property.  
The search of the term that makes the machinery works is one the most challenging issue. In this particular case we will see that one can add
 $${\mathfrak {R}}\circ((-\Delta)^{1/4}Q)-({\mathfrak {R}}\circ((-\Delta)^{1/4}Q)^{\ast}={\mathfrak {R}}\circ((-\Delta)^{1/4}Q)+(-\Delta)^{1/4}Q\circ {\mathfrak {R}}\quad.$$

 We observe  that  $$S_Q(x,y)+\frac{1}{\pi}\frac{1}{x-y}\left[\int_{\R}\frac{Q(y)-Q(z)}{|y-z|^{3/2}}+\int_{\R}\frac{Q(x)-Q(z)}{|x-z|^{3/2}}\right]$$ is the  Schwarz-Kernel of
 $T_R^Q-(T_R^Q)^*$.
 Actually we have
 we have also
\be
\label{II.6-b1}
\begin{array}{l}
\ds T^Q_R v(x)=\frac{1}{\pi}\int_{\R}\frac{(d^{1/2}Q)v(x)-(d^{1/2}Q)v(y)}{x-y}\ dy-{\mathfrak R}\lf[\int_{\R}\frac{Q(\cdot)-Q(z)}{|\cdot-z|^{3/2}}\ dz\ v(\cdot)\rg]\\[5mm]
\ds \quad=\frac{1}{\pi}\int_{\R}\frac{(d^{1/2}Q)v(x)-(d^{1/2}Q)v(y)}{x-y}\ dy-\frac{1}{\pi}\int_{\R}\frac{dy}{x-y}\int_{\R}\frac{Q(x)-Q(z)}{|x-z|^{3/2}}\ dz\ v(x)\\[5mm]
\ds\quad+\frac{1}{\pi}\int_{\R}\frac{1}{x-y}\int_{\R}\frac{Q(y)-Q(z)}{|y-z|^{3/2}}\ dz\ v(y)\ dy\quad.
\end{array}
\ee
Hence
\be
\label{II.6-b2}
\begin{array}{l}
\ds K_{T^Q_R}(x,y)=\frac{1}{\pi}\int_{\R}\frac{dz}{x-z}\lf[\frac{Q(y)-Q(x)}{|y-x|^{3/2}}-\frac{Q(y)-Q(z)}{|z-y|^{3/2}}\rg]+\frac{1}{\pi}\frac{1}{x-y}\int_{\R}\frac{Q(y)-Q(z)}{|y-z|^{3/2}}\ dz\quad.
\end{array}
\ee
In particular
\be
\label{II.6-b3}
\begin{array}{l}
\ds \int_{{\R}}K_{T^Q_R}(x,y)\ dy=-\frac{1}{\pi}\int_{\R}\frac{dy}{x-y}\int_{\R} \frac{Q(z)-Q(y)}{|z-y|^{3/2}}\ dz+\frac{1}{\pi}\int_{{\R}}\frac{dy}{x-y}\int_{\R}\frac{Q(y)-Q(z)}{|y-z|^{3/2}}\ dz\\[5mm]
\quad\quad\quad\quad\ds=-\, 2\, {\mathfrak R}((-\Delta)^{1/4}Q)(x)\quad,
\end{array}
\ee
and
 \be
\label{II.6-b4}
\begin{array}{l}
\ds \int_{{\R}}K_{T^Q_R}(x,y)\ dx=0\quad.
\end{array}
\ee
\medskip

 \subsection{Preliminary estimates of the Kernel of $R^Q-L^Q$}
 In this section we estimate the kernel $$S_Q(x,y)+\frac{1}{\pi}\frac{1}{x-y}\left[\int_{\R}\frac{Q(y)-Q(z)}{|y-z|^{3/2}}+\int_{\R}\frac{Q(x)-Q(z)}{|x-z|^{3/2}}\right]$$ which is exactly 
 the kernel $K_{T_R^Q}(x,y)-K_{T_R^Q}(y,x)$ of the operator $T_R^Q-(T_R^Q)^*,$ in particular we are going to show that it belongs to the functional space
$A_p^{-1/2+1/p}(\R^{2})$ for every $p\ge 2.$\par

\noindent{\bf Step 1. {Estimate of $K_{T_R^Q}(x,y)$}.}\par

  To that aim we decompose $$K_{T_R^Q}(x,y)=K^-_{T_R^Q}(x,y)+K^1_{T_R^Q}(x,y)+K^+_{T_R^Q}(x,y)$$ where
\begin{eqnarray*}
K^-_{T_R^Q}(x,y)&:=&\frac{1}{\pi}\int_{2\,|x-z|<|x-y|}\frac{dz}{x-z}\lf[\frac{Q(y)-Q(x)}{|y-x|^{3/2}}-\frac{Q(y)-Q(z)}{|z-y|^{3/2}}\rg]\\&&+\frac{1}{\pi}\frac{1}{x-y}\int_{2\,|x-z|<|x-y|}\frac{Q(y)-Q(z)}{|y-z|^{3/2}}\ dz
 \\
K^+_{T_R^Q}(x,y)&:=&\frac{1}{\pi}\int_{2\,|x-y|<|x-z|}\frac{dz}{x-z}\lf[\frac{Q(y)-Q(x)}{|y-x|^{3/2}}-\frac{Q(y)-Q(z)}{|z-y|^{3/2}}\rg]\\&&+\frac{1}{\pi}\frac{1}{x-y}\int_{2\,|x-y|<|x-z|}\frac{Q(y)-Q(z)}{|y-z|^{3/2}}\ dz
\end{eqnarray*} 
{\bf {Estimate of $\pi K^-_{T_R^Q}(x,y)$}}\par
We have
\be
\label{II.7}
\begin{array}{l}
\ds \pi K^-_{T_R^Q}(x,y)=\int_{2\,|x-z|<|x-y|}\frac{Q(y)-Q(x)}{x-z}\lf[\frac{1}{|y-x|^{3/2}}-\frac{1}{|z-y|^{3/2}}\rg]\ dz\\[5mm]
\ds\quad\quad+\int_{2\,|x-z|<|x-y|}\frac{Q(z)-Q(x)}{x-z}\ \frac{dz}{|z-y|^{3/2}}+\frac{1}{x-y}\int_{2\,|x-z|<|x-y|}\frac{Q(y)-Q(z)}{|y-z|^{3/2}}\ dz\quad.
\end{array}
\ee
Hence we have
\be
\label{II.8}
\begin{array}{l}
\ds \lf|\pi K^-_{T_R^Q}(x,y)-\int_{2\,|x-z|<|x-y|}\frac{Q(z)-Q(x)}{x-z}\ \frac{dz}{|x-y|^{3/2}}\rg|\\[5mm]
\ds\quad\quad\le|Q(y)-Q(x)|\ \int_{2\,|x-z|<|x-y|}\frac{1}{|x-z|}\frac{1}{|y-x|^{3/2}}\ \frac{|x-z|}{|x-y|}\ dz\\[5mm]
\ds \quad\quad+\frac{1}{|x-y|^{3/2}} \int_{2\,|x-z|<|x-y|}\frac{|Q(z)-Q(x)|}{|x-z|}\ \lf| 1-\frac{1}{\lf|1+\frac{z-x}{x-y}\rg|^{3/2}}\rg|\ dz\\[5mm]
\ds \quad\quad+\frac{1}{|x-y|^{5/2}}\int_{2\,|x-z|<|x-y|}|Q(y)-Q(z)|\ dz\quad.

\end{array}
\ee
Using the inequality of $\lf| 1-\frac{1}{\lf|1+\frac{z-x}{x-y}\rg|^{3/2}}\rg|\le C\,\lf|\frac{z-x}{x-y}\rg|$ in the range $|x-y|>2\, |x-z|$ we have
\be
\label{II.9}
\begin{array}{l}
\ds \lf|\pi K^-_{T_R^Q}(x,y)-\int_{2\,|x-z|<|x-y|}\frac{Q(z)-Q(x)}{x-z}\ \frac{dz}{|x-y|^{3/2}}\rg|\lesssim\, \left[\ \frac{|Q(y)-Q(x)|}{|x-y|^{3/2}}\right.\\[5mm]
\ds\quad\quad+\left.\frac{1}{|x-y|^{5/2}}\, \int_{2\,|x-z|<|x-y|}[|Q(x)-Q(z)|+|Q(y)-Q(z)|]\ dz\right]
\end{array}
\ee
{\bf {Estimate of  $\pi K^+_{T_R^Q}(x,y)$.}} We use the fact that
\be
\label{II.11}
 \int_{2\,|x-y|<|x-z|}\frac{dz}{x-z}\frac{Q(y)-Q(x)}{|y-x|^{3/2}}=0\quad,
\ee
and $$\{2\,|x-y|<|x-z|\}\subseteq \{|x-y|<|y-z|\}\quad.$$
We write
\begin{eqnarray}
\label{II.12}
\ds \lf|\pi K^+_{T_R^Q}(x,y)\rg|&=&\lf| \int_{2\,|x-y|<|x-z|}\frac{dz}{x-z} \frac{Q(y)-Q(z)}{|z-y|^{3/2}} +\frac{1}{x-y}
\int_{2\,|x-y|<|x-z|}  \frac{Q(y)-Q(z)}{|z-y|^{3/2}} dz\rg|\nonumber\\[5mm]
&\lesssim& 
\int_{|x-y|<|y-z|}\frac{dz}{|x-z|} \frac{|Q(y)-Q(z)|}{|z-y|^{3/2}} +\frac{1}{x-y}
\int_{|x-y|<|y-z|}  \frac{|Q(y)-Q(z)|}{|z-y|^{3/2}} dz  \nonumber\\[5mm]
&\lesssim&
   \int_{|v|>|x-y|}\frac{1}{|v|\left[1+\frac{|x-y|}{|v|}\right]} \frac{|Q(y)-Q(y+v)|}{|v|^{3/2}} dv\\[5mm]
   &+&\frac{1}{|x-y|}\int_{|v|>|x-y|}\frac{|Q(y)-Q(y+v)|}{|v|^{3/2}} dv.\nonumber
 \end{eqnarray}

{\bf Estimate of $\pi K^1_{T_R^Q}(x,y)$}\par 
 We have $\pi K^1_{T_R^Q}(x,y)=\pi K^{1,+}_{T_R^Q}(x,y)+\pi K^{1,-}_{T_R^Q}(x,y)$ where
\begin{eqnarray*}
\pi K^{1,-}_{T_R^Q}(x,y)&=&-\int_{\{2\, |y-z|< |x-z|\}\cap\{2^{-1}\,|x-y|<|x-z|<2\, |x-y|\}}\frac{dz}{x-z}\frac{Q(y)-Q(z)}{|z-y|^{3/2}}\\
&+&\frac{1}{x-y}\int_{\{2\, |y-z|< |x-z|\}\cap\{2^{-1}\,|x-y|<|x-z|<2\, |x-y|\}}\ \frac{Q(y)-Q(z)}{|z-y|^{3/2}}\ dz\quad,
 \end{eqnarray*}
and
\[
\begin{array}{l}
\ds \pi K^{1,+}_{T_R^Q}(x,y)=\int_{\{|x-z|<2\, |y-z|\}\cap\{2^{-1}\,|x-y|<|x-z|<2\, |x-y|\} }\frac{dz}{x-z}\frac{Q(y)-Q(x)}{|y-x|^{3/2}}\\[5mm]
\ds\quad\quad-\int_{\{|x-z|<2\, |y-z|\}\cap\{2^{-1}\,|x-y|<|x-z|<2\, |x-y|\} }\frac{dz}{x-z}\frac{Q(y)-Q(z)}{|z-y|^{3/2}}\\[5mm]
\ds\quad\quad+\frac{1}{x-y}\int_{\{|x-z|<2\, |y-z|\}\cap\{2^{-1}\,|x-y|<|x-z|<2\, |x-y|\} }\ \frac{Q(y)-Q(z)}{|z-y|^{3/2}}\ dz\quad.
\end{array}
\]
We have obviously
\be
\label{II.10}
\lf|\pi K^{1,+}_{T_R^Q}(x,y)\rg|\le C\, \frac{|Q(x)-Q(y)|}{|x-y|^{3/2}}+C\, |x-y|^{-5/2}\ \int_{2^{-1}|x-y|\le |x-z|\le 2\, |x-y|}|Q(y)-Q(z)|\ dz\quad,
\ee
and
\be
\label{II.10-a}
\begin{array}{l}
\ds |\pi K^{1,-}_{T_R^Q}(x,y)|=\lf|\int_{2\, |y-z|< |x-z|}\frac{Q(y)-Q(z)}{|z-y|^{3/2}}\ \lf[-\frac{1}{x-z}+\frac{1}{x-y}\rg]\ dz\rg|\\[5mm]
\ds\quad\quad\le\frac{1}{|x-y|^2} \int_{2\, |y-z|< |x-z|}\lf|\frac{Q(y)-Q(z)}{|z-y|^{1/2}}\rg| \quad.
\end{array}
\ee
{\bf Step 2: Estimate of $\pi (K_{T_R^Q}(x,y)-K_{T_R^Q}(y,x))$}\par
 
  In \eqref{II.9} we have subtracted  $\int_{2\,|x-z|<|x-y|}\frac{Q(z)-Q(x)}{x-z}\ \frac{dz}{|x-y|^{3/2}}$.
  Since we are considering $ \pi (K_{T_R^Q}(x,y)-K_{T_R^Q}(y,x))$ this means that we have to estimate

\begin{eqnarray}\label{II.10-b}
&&\int_{2\,|x-z|<|x-y|}\frac{Q(z)-Q(x)}{x-z}\ \frac{dz}{|x-y|^{3/2}}-\int_{2\,|y-z|<|x-y|}\frac{Q(z)-Q(y)}{y-z}\ \frac{dz}{|x-y|^{3/2}}\nonumber\\
&=&- \frac{{\mathfrak R}Q(x) }{|x-y|^{3/2}} +\frac{{\mathfrak R}Q(y) }{|x-y|^{3/2}}\\
&&-\int_{2\,|x-z|>|x-y|}\frac{Q(z)-Q(x)}{x-z}\ \frac{dz}{|x-y|^{3/2}}+\int_{2\,|y-z|>|x-y|}\frac{Q(z)-Q(y)}{x-z}\ \frac{dz}{|x-y|^{3/2}}\nonumber.
\end{eqnarray}
  Therefore we first estimate 
$$ \left|\underbrace{\int_{2\,|x-z|>|x-y|}\frac{Q(z)-Q(x)}{x-z}\ \frac{dz}{|x-y|^{3/2}}}_{(1)}-\underbrace{\int_{2\,|y-z|>|x-y|}\frac{Q(z)-Q(y)}{y-z}\ \frac{dz}{|x-y|^{3/2}}}_{(2)}\right|$$
We define the following sets:
$$A_1:=\{2^{-1}|x-y|\le |x-z|<2|x-y|\},~~A_2:=\{ 2|x-y| < |x-z|,~|x-y|<|y-z|<2|x-y| \}$$
and $$A_3=B_3:=\{|z-x|>2|x-y|,~|z-y|>2|x-y|\}$$
and
in a similar way
$$B_1:=\{2^{-1}|x-y|\le |y-z|<2|x-y|\}~~B_2:=\{ 2|x-y| < |y-z|,~|x-y|<|x-z|<2|x-y| \}\quad, $$
We split $(1)$ as follows:
\begin{eqnarray}
&&\int_{2\,|x-z|>|x-y|}\frac{Q(z)-Q(x)}{x-z}\ \frac{dz}{|x-y|^{3/2}}=\int_{A_1}\frac{Q(z)-Q(x)}{x-z}\ \frac{dz}{|x-y|^{3/2}}\nonumber\\
&&+ \int_{A_2}\frac{Q(z)-Q(x)}{x-z}\ \frac{dz}{|x-y|^{3/2}}+\int_{A_3}\frac{Q(z)-Q(x)}{x-z}\ \frac{dz}{|x-y|^{3/2}}
\end{eqnarray}
The following estimates hold:
\begin{eqnarray}
\left|\int_{A_1}\frac{Q(z)-Q(x)}{x-z}\ \frac{dz}{|x-y|^{3/2}}\right|&\lesssim&  {|x-y|^{-5/2}}\int_{A_1}|Q(z)-Q(x)|dz.\label{A1}\\
\left|\int_{A_2}\frac{Q(z)-Q(x)}{x-z}\ \frac{dz}{|x-y|^{3/2}}\right|&\lesssim & \left|\int_{A_2}\left[\frac{Q(z)-Q(x)}{x-z}-\frac{Q(y)-Q(x)}{x-z}\right]\frac{dz}{|x-y|^{3/2}}\right| \label{A2}\\[5mm]
&\lesssim& \frac{1}{|x-y|^{3/2}}\int_{2|x-y|>|y-z|}\frac{|Q(z)-Q(y)|}{|y-z|}dz.\nonumber\end{eqnarray}
In an analogous way we find for $(2)$ the following:
\begin{equation}\label{B1}
\left|\int_{B_1}\frac{Q(z)-Q(y)}{y-z}\ \frac{dz}{|x-y|^{3/2}}\right|\lesssim  {|x-y|^{-5/2}}\int_{B_1}|Q(z)-Q(y)|dz.
\end{equation}
\begin{eqnarray}\label{B2}
\left|\int_{B_2}\frac{Q(z)-Q(y)}{x-z}\ \frac{dz}{|x-y|^{3/2}}\right|&\lesssim & \left|\int_{B_2}\left[\frac{Q(z)-Q(y)}{y-z}-\frac{Q(x)-Q(y)}{y-z}\right]\frac{dz}{|x-y|^{3/2}}\right|\nonumber\\[5mm]
&\lesssim& \frac{1}{|x-y|^{3/2}}\int_{2|x-y|>|x-z|}\frac{|Q(z)-Q(x)|}{|x-z|}dz.\end{eqnarray}
In order to estimate $\int_{A_3}\frac{Q(z)-Q(x)}{x-z}\ \frac{dz}{|x-y|^{3/2}}$ we need to put it together with
$\int_{A_3}\frac{Q(z)-Q(y)}{y-z}\ \frac{dz}{|x-y|^{3/2}}$, which is a sort of {\bf compensation effect}.
\par
We have
\begin{eqnarray}\label{A3}
&&\left|\int_{A_3}\frac{Q(z)-Q(x)}{x-z}\ \frac{dz}{|x-y|^{3/2}}-\int_{A_3}\frac{Q(z)-Q(y)}{y-z}\ \frac{dz}{|x-y|^{3/2}}\right|\nonumber\\[5mm]
&&=\left|\int_{A_3}Q(z)-Q(x)\left[\frac{1}{x-z}-\frac{1}{y-z}\right]\ \frac{dz}{|x-y|^{3/2}}\right|\\[5mm]
&&\lesssim\frac{1}{|x-y|^{5/2}}\int_{|x-z|>2|x-y|,|y-z|>2|x-y|}|Q(z)-Q(x)|dz\quad.\nonumber
\end{eqnarray}
Combining (\ref{II.7})-(\ref{A3}) gives
\be
\label{II.13}
\begin{array}{l}
\ds \lf|\pi (K_{T_R^Q}(x,y)-K_{T_R^Q}(y,x))\rg|\lesssim  \frac{|Q(x)-Q(y)|}{|x-y|^{3/2}}+ \frac{|{\mathfrak R}Q(x)-{\mathfrak R}Q(y)|}{|x-y|^{3/2}}\\[5mm]
\ds \quad\quad+ \left|\int_{2\,|x-z|>|x-y|}\frac{Q(z)-Q(x)}{x-z}\ \frac{dz}{|x-y|^{3/2}}-\int_{2\,|x-z|>|x-y|}\frac{Q(z)-Q(x)}{x-z}\ \frac{dz}{|x-y|^{3/2}}\right| \\[5mm]
\ds \quad\quad+\, |x-y|^{-5/2}\ \int_{2^{-1}|x-y|\le |x-z|\le 2\, |x-y|}|Q(y)-Q(z)|\ dz\\[5mm]
\ds \quad\quad+\, |x-y|^{-5/2}\ \int_{2^{-1}|x-y|\le |y-z|\le 2\, |x-y|}|Q(x)-Q(z)|\ dz\\[5mm]
\ds \quad\quad+\frac{1}{|x-y|^{5/2}}\, \int_{2\,|x-z|<|x-y|}|Q(x)-Q(z)|\ dz+\frac{1}{|x-y|^{5/2}}\, \int_{2\,|y-z|<|x-y|}|Q(y)-Q(z)|\ dz\\[5mm]
\ds\quad\quad+\frac{1}{|x-y|^2} \lf[\int_{2\, |y-z|< |x-z|}\lf|\frac{Q(y)-Q(z)}{|z-y|^{1/2}}\rg| \ dz+\int_{2\, |x-z|< |y-z|}\lf|\frac{Q(x)-Q(z)}{|z-x|^{1/2}}\rg|\ dz \rg]\\[5mm]
\ds\quad\quad+\int_{|x-z|\ge2 |x-y|}\lf|\frac{Q(z)-Q(x)}{x-z}-\frac{Q(z)-Q(y)}{y-z}\rg|\ \frac{dz}{|x-y|^{3/2}}\\[5mm]
\ds\quad\quad + \int_{|v|>|x-y|}\frac{1}{|v|\left[1+\frac{|x-y|}{|v|}\right]} \left[\frac{|Q(y)-Q(y+v)|}{|v|^{3/2}} +\frac{|Q(x)-Q(x+v)|}{|v|^{3/2}}\right] dv\\[5mm]
   \ds\quad\quad +\frac{1}{|x-y|}\int_{|v|>|x-y|}\left[\frac{|Q(y)-Q(y+v)|}{|v|^{3/2}}+\frac{|Q(x)-Q(x+v)|}{|v|^{3/2}}\right] dv\quad.\
\end{array}
\ee
 Using (\ref{II.13}) we shall now prove  lemma~\ref{lm-II.1} which is the goal of the present subsection.
 We aim at proving that  
 \begin{equation}\label{II.13b}
 {\mathfrak{S}}^Q(x,y):=K_{T_R^Q}(x,y)-K_{T_R^Q}(y,x)  \end{equation}
satisfies
for any $p\ge 2$ the following bound
\be
\label{II.14}
\|{\mathfrak{S}}^Q(x,y)\|_{A_{p,2}^{-1/2+1/p}}=\lf[\int_{\R} h^{2-2/p}\ \|{\mathfrak{S}}^Q(\cdot,\cdot+h)\|^2_p\ dh\rg]^{1/2}\le C_p\ \|Q\|_{\dot{H}^{1/2}({\R})}\quad.
\ee
 
The estimates \eqref{II.14} will imply that ${\mathcal T}_{{\mathfrak{S}}^Q}$ defines an {\bf abstract multi-commutator}  satisfying the compensation property of lemma~\ref{lm-I.1}. It is explicitly given by\footnote{We are using (\ref{II.6-b3}) and (\ref{II.6-b4}). In particular
we have that
$$\int_{\R} {\mathfrak{S}}^Q(y,x) dy= <T_R^Q-(T_R^Q)^\ast(x),1>=-2{\mathfrak R}\circ(-\Delta)^{1/4}Q.$$} 
\begin{eqnarray*}
\label{II.14-a}
\ds{\mathcal T}_{{\mathfrak{S}}^Q}&=&{\mathfrak R}\circ(Q\circ(-\Delta)^{1/4}-(-\Delta)^{1/4}\circ Q-(-\Delta)^{1/4}Q)\\[5mm]&&-(Q\circ(-\Delta)^{1/4}-(-\Delta)^{1/4}\circ Q-(-\Delta)^{1/4}Q)\circ{\mathfrak R}
\\[5mm]&&-2\,(-\Delta)^{1/4}Q\circ{\mathfrak R}-2\,{\mathfrak R}\circ((-\Delta)^{1/4}Q).
\end{eqnarray*}
 
 \noindent{\bf Proof of Lemma~\ref{lm-II.1}.} We bound the $A_{p,2}^{-1/2+1/p}$ norm of each terms in the r.h.s. of (\ref{II.13}).
We first have
\be
\label{II.15}
\begin{array}{l}
\ds\int_{\R} h^{2-2/p}\ \lf[\int_{\R} \frac{|Q(x)-Q(x+h)|^p}{|h|^{3p/2}}dx \rg]^{2/p}\ dh=\int_{\R}\frac{dh}{h^{1+2/p}}\ \|Q(\cdot)-Q(\cdot+h)\|^2_{L^p}\\[5mm]
\ds\quad\quad=\|Q\|^2_{B^{1/p}_{p,2}(\R)}\le C_p\ \|Q\|^2_{\dot{H}^{1/2}({\R})}\quad.
\end{array}
\ee
We have also
\be
\label{II.16}
\begin{array}{l}
\ds\int_{\R} h^{2-2/p}\ \lf[\int_{\R} \frac{|{\mathfrak R}Q(x)-{\mathfrak R}Q(x+h)|^p}{|h|^{3p/2}} \rg]^{2/p}\ dh=\int_{\R}\frac{dh}{h^{1+2/p}}\ \|{\mathfrak R}Q(\cdot)-{\mathfrak R}Q(\cdot+h)\|^2_{L^p}\\[5mm]
\ds\quad\quad=\|{\mathfrak R}Q\|^2_{B^{1/p}_{p,2}(\R)}\le C_p\ \|Q\|^2_{B^{1/p}_{p,2}(\R)}\le C_p\ \|Q\|^2_{\dot{H}^{1/2}({\R})}\quad.
\end{array}
\ee
We have
\be
\label{II.17}
\begin{array}{l}
\ds\int_{\R} h^{2-2/p}\ \lf[\int_{\R} h^{-5p/2}\lf|\int_{2^{-1}|h|\le |x-z|\le 2\, |h|}|Q(x)-Q(z)|\ dz\rg|^p\ dx\rg]^{2/p}\ dh\\[5mm]
\ds\quad\quad\le \int_{\R} h^{-3-2/p}\ \lf\|\int_{1/2}^2|Q(\cdot)-Q(\cdot+t\,h)|\ |h|\ dt \rg\|^2_{L^p}\ dh\\[5mm]
\ds\quad\quad\le  \int_{\R} h^{-1-2/p}\ \lf[\int_{1/2}^2\lf\|Q(\cdot)-Q(\cdot+t\,h)\rg\|_{L^p} dt\rg]^2\ dh\\[5mm]
\ds\quad\quad\le  2\ \int_{1/2}^2dt\int_{\R} h^{-1-2/p}\ \lf\|Q(\cdot)-Q(\cdot+t\,h)\rg\|_{L^p}^2\ dh\\[5mm]
\ds\quad\quad\le C\ \int_{\R} h^{-1-2/p}\ \lf\|Q(\cdot)-Q(\cdot+h)\rg\|_{L^p}^2\ dh\le C_p\ \|Q\|^2_{B^{1/p}_{p,2}(\R)}\le C_p\ \|Q\|^2_{\dot{H}^{1/2}({\R})}\quad,
\end{array}
\ee
where we have used successively {\it Minkowski integral inequality} and {\it Cauchy Schwartz inequality}. 

We have
\be
\label{II.17-a}
\begin{array}{l}
\ds\int_{\R} h^{2-2/p}\ \lf[\int_{\R} |h|^{-2p} \lf|\int_{2\, |x+h-z|< |x-z|}\lf|\frac{Q(x+h)-Q(z)}{|x+h-z|^{1/2}}\rg| dz\rg|^p\ dx\rg]^{2/p}\ dh\\[5mm]
\ds\le\int_{\R} h^{2-2/p}\ \lf[\int_{\R} |h|^{-2p} \lf|\int_{\, |th|< |h|}\lf|\frac{Q(x+h)-Q(x+h+th)}{|th|^{1/2}}\rg| d(th)\rg|^p\ dx\rg]^{2/p}\ dh\\[5mm]
\ds\le\int_{\R} h^{-1-2/p}\ \lf[\int_{0}^{1} \frac{dt}{\sqrt{t}}\|Q(\cdot)-Q(\cdot+th)\|_{L^p}\rg]^2\ dh\\[5mm]
\ds\le \int_0^1\frac{dt}{t^{1/2-2/p}}\int_{\R} (th)^{-1-2/p}\  \|Q(\cdot)-Q(\cdot+th)\|^2_{L^p}\ d(th)\le C_p\ \|Q\|^2_{B^{1/p}_{p,2}(\R)}\le C_p\ \|Q\|^2_{\dot{H}^{1/2}({\R})}\ .
\end{array}
\ee
We have also
\be
\label{II.18}
\begin{array}{l}
\ds\int_{\R} h^{2-2/p}\ \lf\|h^{-5/2}\int_{2\,|x-z|<|h|}|Q(x)-Q(z)| \ dz\rg\|^2_{L^p}\ dh\\[5mm]
\ds\quad\quad\le \int_{\R} h^{-3-2/p}\ \lf\|\int^{1/2}_0|Q(\cdot)-Q(\cdot+t\,h)|\ |h|\ dt \rg\|^2_{L^p}\ dh\\[5mm]
\ds\quad\quad\le  \int_{\R} h^{-1-2/p}\ \lf[\int^{1/2}_0\lf\|Q(\cdot)-Q(\cdot+t\,h)\rg\|_{L^p} dt\rg]^2\ dh\\[5mm]
\ds\quad\quad\le  \int^{1/2}_0dt\int_{\R} h^{-1-2/p}\ \lf\|Q(\cdot)-Q(\cdot+t\,h)\rg\|_{L^p}^2\ dh\\[5mm]
\ds\quad\quad\le  \int^{1/2}_0t^{2/p}dt\int_{\R} (th)^{-1-2/p}\ \lf\|Q(\cdot)-Q(\cdot+t\,h)\rg\|_{L^p}^2\ d(th)\\[5mm]
\ds\quad\quad\le C\ \int_{\R} h^{-1-2/p}\ \lf\|Q(\cdot)-Q(\cdot+h)\rg\|_{L^p}^2\ dh\le C\ \|Q\|^2_{B^{1/p}_{p,2}(\R)}\le C\ \|Q\|^2_{\dot{H}^{1/2}({\R})}\quad.
\end{array}
\ee
 Using the fact that
\be
\label{II.18-a}
\int_{2\,|y-z|\ge |h|}\frac{Q(y)-Q(x)}{y-z}\ dz=0
\ee
we have
\be
\label{II.19}
\begin{array}{l}
\ds\int_{\R} h^{2-2/p}\ \lf\|h^{-3/2}\int_{|x-z|\ge 2\,|h|}\lf|\frac{Q(z)-Q(x)}{x-z}-\frac{Q(z)-Q(x+h)}{x+h-z}\rg|\ {dz}\rg\|^2_{L^p}\ dh\\[5mm]
\ds\quad\quad\le\int_{\R} h^{-1-2/p}\ \lf\|\int_{|x-z|\ge 2\, |h|}\lf|{Q(z)-Q(x)}\ \lf[ \frac{1}{x-z}-\frac{1}{x+h-z}  \rg]\rg|\ {dz}\rg\|^2_{L^p}\ dh\\[5mm]
\ds\quad\quad+\int_{\R} h^{2-2/p}\ \lf[\int_{\R} h^{-5p/2}\lf|\int_{2^{-1}|h|\le |x-z|\le 2\, |h|}|Q(x)-Q(z)|\ dz\rg|^p\ dx\rg]^{2/p}\ dh
\end{array}
\ee
The second term of the right-hand side of (\ref{II.19}) has already be controlled in (\ref{II.17}). Hence we bound now
\be
\label{II.20}
\begin{array}{l}
\ds\int_{\R} h^{-1-2/p}\ \lf\|\int_{|x-z|\ge 2\, |h|}\lf|{Q(z)-Q(x)}\ \lf[ \frac{1}{x-z}-\frac{1}{x+h-z}  \rg]\rg|\ {dz}\rg\|^2_{L^p}\ dh\\[5mm]
\ds\quad\quad\le \int_{\R} h^{-1-2/p} \lf\|\int_{|x-z|\ge 2\, |h|}\lf|{Q(z)-Q(x)}\rg|\ \frac{|h|}{|x-z|^2}{dz}\rg\|^2_{L^p}\ dh\\[5mm]
\ds\quad\quad\le \int_{\R} h^{-1-2/p} \lf\|\int_{t\,|h|\ge 2\, |h|}\lf|{Q(x+th)-Q(x)}\rg|\ \frac{|h|}{|th|^2}{d(th)}\rg\|^2_{L^p}\ dh\\[5mm]
\ds\quad\quad\le \int_{\R} h^{-1-2/p}\ \lf[\int_{2}^{+\infty}\ \frac{dt}{t^2}\lf\|Q(x+th)-Q(x)\rg\|_{L^p}\rg]^2\ dh\\[5mm]
\ds\quad\quad\le \int_{\R} h^{-1-2/p}\ \int_{2}^{+\infty}\ \frac{dt}{t^{5/2}}\ \lf\|Q(x+th)-Q(x)\rg\|_{L^p}^2\ dh\ \int_{2}^{+\infty}\ \frac{dt}{t^{3/2}}\\[5mm]
\ds\quad\quad\le C\, \int_{2}^{+\infty}\ \frac{dt}{t^{5/2-2/p}}\int_{\R} (th)^{-1-2/p}\   \lf\|Q(x+th)-Q(x)\rg\|_{L^p}^2\ d(th)\\[5mm]
\ds\quad\quad\le C\ \|Q\|^2_{B^{1/p}_{p,2}(\R)}\le C\ \|Q\|^2_{\dot{H}^{1/2}({\R})}\quad.
\end{array}
\ee
Combining (\ref{II.19}), (\ref{II.17}) and (\ref{II.20}) we finally obtain
\be
\label{II.21}
\ds\int_{\R} h^{2-2/p}\ \lf\|h^{-3/2}\int_{|x-z|\ge 2\,|h|}\lf|\frac{Q(z)-Q(x)}{x-z}-\frac{Q(z)-Q(x+h)}{x+h-z}\rg|\ {dz}\rg\|^2_{L^p}\ dh\le C\ \|Q\|^2_{\dot{H}^{1/2}({\R})}\quad.
\ee
To conclude the proof of lemma~\ref{lm-II.1} we have to bound
  \be
\label{II.22}
\begin{array}{l}
 \ds\int_{\R} h^{2-2/p}\left( \int_{\R}\left| \int_{|v|>|h|}\frac{1}{|v|\left[1+\frac{|h|}{|v|}\right]} \frac{|Q(y)-Q(y+v)|}{|v|^{3/2}}dv\right|^p dx\right)^{2/p} dh\\
 + \ds\int_{\R} h^{2-2/p}\left( \int_{\R}\left|  \int_{|v|>|h|}\frac{1}{|v|\left[1+\frac{|h|}{|v|}\right]} \frac{|Q(x)-Q(x+v)|}{|v|^{3/2}} dv\right|^p dx\right)^{2/p} dh\quad.
\end{array}
\ee
and
 \be
\label{II.23}
\begin{array}{l}
 \ds\int_{\R} h^{2-2/p}\left( \int_{\R}\left|\int_{|v|>|h|}\frac{1}{|h|}   \frac{|Q(y)-Q(y+v)|}{|v|^{3/2}}  dv\right|^p dx\right)^{2/p} dh\\
 + \ds\int_{\R} h^{2-2/p}\left( \int_{\R}\left| \int_{|v|>|h|}\frac{1}{|h|} \frac{|Q(x)-Q(x+v)|}{|v|^{3/2}}  dv\right|^p dx\right)^{2/p} dh\quad.
 \end{array}
\ee
We are going to estimate only the first term in \eqref{II.22} (the other terms can be estimated in a similar way). Let fix $0<\varepsilon< 1/2-1/p$.
 \be
\label{II.24}
\begin{array}{l}
\ds\int_{\R} h^{2-2/p}\left( \int_{\R}\left| \int_{|v|>|h|}\frac{1}{|v|\left[1+\frac{|h|}{|v|}\right]} \frac{|Q(y)-Q(y+v)|}{|v|^{3/2}}dv\right|^p dx\right)^{2/p} dh\\[5mm]
 \le \ds\int_{\R} h^{2-2/p}\left( \int_{\R} \left| \int_{t>1}\frac{1}{|th|\left[1+\frac{|h|}{|th|}\right]} \frac{|Q(x+h)-Q(x+(t+1)h)|}{|th|^{3/2}}|h|dt \right|^p dx\right)^{2/p} dh\\[5mm]
 = \ds\int_{\R} h^{-1-2/p}\left( \int_{\R} \left| \int_{t>1}\frac{1}{|t|\left[1+\frac{1}{|t|}\right]} \frac{|Q(x')-Q(x'+th)|}{|t|^{3/2}} dt \right|^p dx'\right)^{2/p} dh\\[5mm]
 \lesssim \ds\int_{\R} h^{-1-2/p}\left( \int_{t>1}t^{-3/2}\|Q(\cdot)-Q(\cdot+th)\|_{L^p}\right)^{2} dh\\[5mm]
 \lesssim \ds\int_{\R} h^{-1-2/p}\left(\int_{t>1}t^{-1-2\varepsilon} dt\right)\int_{t>1}t^{-2+2\varepsilon}\|Q(\cdot)-Q(\cdot+th)\|_{L^p}^2 dt\\[5mm]
 \lesssim (\int_{t>1}t^{-2+2\varepsilon}\left(\int_{\R}( th)^{-1-2/p}t^{1+2/p}t^{-1}\|Q(\cdot )-Q(\cdot+th)\|_{L^p}^2 d(ht) \right)dt\\[5mm]
 =(\int_{t>1}t^{-2+2\varepsilon+2/p}dt \left(\int_{\R}( th)^{-1-2/p} \|Q(\cdot )-Q(\cdot+th)\|_{L^p}^2 d(ht) \right)\\[5mm]
\lesssim C_p\ \|Q\|^2_{B^{1/p}_{p,2}(\R)}\le C_p\ \|Q\|^2_{\dot{H}^{1/2}({\R})}\quad .
 \end{array}
\ee
We observe that the integral $\int_{t>1}t^{-2+2\varepsilon-2/p}dt$ in \eqref{II.24} converges since  $\varepsilon< 1/2-1/p$.
  
 Combining (\ref{II.13}), (\ref{II.15}), (\ref{II.16}), (\ref{II.17}), (\ref{II.17-a}), (\ref{II.18}), (\ref{II.21}) and (\ref{II.22})-(\ref{II.24}) we obtain (\ref{II.14}) and lemma~\ref{lm-II.1} is proved.\hfill $\Box$

\medskip

\begin{Rm}
\label{rem-sect-1} We already knew that ${\mathcal T}_{{\mathfrak{S}}^Q}$ is a multi-commutator since it is given by   
$${\mathfrak R}\circ  {\mathcal{T}}_{K_{d^{1/2} Q}}- {\mathcal{T}}_{K_{d^{1/2} Q}}\circ{\mathfrak R}-\underbrace{2\,(-\Delta)^{1/4}Q\circ{\mathfrak R}-2\,{\mathfrak R}\circ((-\Delta)^{1/4}Q)}_{(1)}$$ and $(1)$ maps $L^2$ into the Hardy space ${{\mathcal H}^1({\R})}$
since for $f,g\in L^2(\R)$ one has the Coifman-Rochberg-Weiss commutator
\[
\|f\,{\mathfrak R}(g)+{\mathfrak R}(f)\, g\|_{{\mathcal H}^1({\R})}\le C\ \|f\|_{L^2}\ \|g\|_{L^2}.
\]
Nevertheless, the information provided by (\ref{II.14})and the $A^s_{p,2}$ bound is new and in particular, thanks to lemma~\ref{stabKernel} , it permits to generate new multi-commutators. Indeed, for any function $P(x)\in L^\infty({\R},\mbox{Sym}_m)$ we still have obviously 
\[
\| P(x)\ {\mathfrak{S}}^Q(x,y)\|_{A_p^{-1/2+1/p}} \le C_p\ \|P\|_\infty \|Q\|_{\dot{H}^{1/2}({\R})}\quad.
\]
If then one considers
\[
W(x,y):=P(x)\ {\mathfrak{S}}^Q(x,y)-(P(y) {\mathfrak{S}}^Q(y,x))^t=P(x)\ {\mathfrak{S}}^Q(x,y)+ {\mathfrak{S}}^Q(x,y)\ P(y)
\]
this generates obviously a new multi-commutator.\hfill $\Box$
\end{Rm}
\medskip

As a matter of illustration of the previous remark, starting from
$$
P(x)K_{ d^{1/2}Q}(x,y)=P(x)\frac{Q(y)-Q(x)}{|x-y|^{3/2}},
 $$
 one considers 
$$P(x)\frac{Q(y)-Q(x)}{|x-y|^{3/2}}+ \frac{Q(y)-Q(x)}{|x-y|^{3/2}}P(y)$$
which is the Schwarz kernel of
\[
P\circ d^{1/2}Q+d^{1/2}Q\circ P
\]
We compute $\int_{\R}P(x)\frac{Q(y)-Q(x)}{|x-y|^{3/2}}+ \frac{Q(y)-Q(x)}{|x-y|^{3/2}}P(y)$ which is given by
\[
\begin{array}{l}
\ds\int_{\R} P(x) \ \frac{Q(y)-Q(x)}{|x-y|^{3/2}} \ dy+\int_{\R}  \frac{Q(y)-Q(x)}{|x-y|^{3/2}}\ P(y) \ dy\\[5mm]
\ds\quad=-P(x)\ (-\Delta)^{1/4}Q+\int_{\R}  \frac{Q(y)\,P(y)-Q(x)\, P(x)+ Q(x)\, P(x)- Q(x)P(y)}{|x-y|^{3/2}} \ dy\\[5mm]
\ds\quad=-P\ (-\Delta)^{1/4}Q- (-\Delta)^{1/4}(QP)+ Q\, (-\Delta)^{1/4}P
\end{array}
\]
Then  we deduce the following lemma
\begin{Lm}
\label{lm-multi} Let $Q\in \dot{H}^{1/2}({\R},\mbox{Sym}_m)$ and $P(x)\in L^\infty({\R},\mbox{Sym}_m)$ then the following operator
\begin{eqnarray*}
V_{P,Q}&:=&(PQ)\circ(-\Delta)^{1/4} - P\circ(-\Delta)^{1/4}\circ Q\\
&+&Q\circ(-\Delta)^{1/4}\circ P -(-\Delta)^{1/4}\circ (QP)-P\ (-\Delta)^{1/4}Q- (-\Delta)^{1/4}(QP)+ Q\, (-\Delta)^{1/4}P\quad.
\end{eqnarray*}
is mapping continuously $L^2$ into $B^{-1/2+1/p}_{2p/(p+2),2}$ for any $2<p$. Since $B^{-1/2+1/p}_{2p/(p+2),2}\hookrightarrow H^{-1/2}$ we have in particular for any $v\in L^2({\R})$
\be
\label{lm-multi-1}
\|V_{P,Q}(v)\|_{H^{-1/2}({\R})}\le C\ \|P\|_{L^\infty({\R})}\ \|Q\|_{\dot{H}^{1/2}({\R})}\ \|v\|_{L^2({\R})}\quad.
\ee
\end{Lm}

\medskip

\subsection{Generating Multi-Commutators from $d^{1/2}{\mathfrak R}_Q:={\mathfrak R}\circ d^{1/2}Q\circ{\mathfrak R}$.}
In this section we are going to generate a multi-commutator starting from 
Let $Q\in \dot{H}^{1/2}({\R})$ we consider 
\[
{\mathfrak R}_{Q}:={\mathfrak R}\circ Q\circ{\mathfrak R}
\]
We have
\be
\label{II-f-1}
{\mathfrak R}_{Q}(v)(x)=\frac{1}{\pi^2}\int_{\R}\int_{\R}\ \frac{Q(z)}{(x-z)\, (z-y)}\ dz\ v(y)\ dy
\ee
We shall denote for $x\ne y$
\be
\label{II-f-2}
{\mathcal R}_Q(x,y):=\frac{1}{\pi^2}\int_{\R}\ \frac{Q(z)}{(x-z)\, (z-y)}\ dz\quad,
\ee
where we observe first that 
\be
\label{II-f-2a}
{\mathcal R}_Q(x,y)={\mathcal R}_Q(y,x)\quad,
\ee
and, for $x<y$,
\be
\label{II-f-3}
\begin{array}{l}
\ds 0=\lim_{\ep\rightarrow 0}\int_{{\R}\setminus B_\ep(x)\cup B_\ep(y)}\ \frac{dz}{(x-z)\, (z-y)}=2\ \lim_{\ep\rightarrow 0}\int_{(-\infty,(x+y)/2)\setminus B_\ep(x)}\ \frac{dz}{(x-z)\, (z-y)}\quad.
\end{array}
\ee
We have indeed
\be
\label{II-f-4}
\begin{array}{l}
\displaystyle\lim_{\ep\to 0}\ds\int_{(-\ep^{-1},(x+y)/2)\setminus B_\ep(x)}\ \frac{dz}{(x-z)\, (z-y)}= \frac{1}{x-y}\lim_{\ep\to 0}\int_{(-\ep^{-1},(x+y)/2)\setminus B_\ep(x)}\ \lf[\frac{1}{x-z} +\frac{1}{z-y}\rg]\ dz=0\quad,\\[5mm]
\end{array}
\ee
and hence, the singular integral (\ref{II-f-2}) has to be understood in the following sense for $x<y$
\be
\label{II-f-5}
{\mathcal R}_Q(x,y):=\frac{1}{\pi^2}\int_{z<(x+y)/2}\ \frac{Q(z)-Q(x)}{(x-z)\, (z-y)}\ dz+\frac{1}{\pi^2}\int_{z>(x+y)/2}\ \frac{Q(z)-Q(y)}{(x-z)\, (z-y)}\ dz\quad.
\ee
We shall now compute and estimate the Schwartz Kernel associated to
\[ {\mathfrak R}_{d^{1/2}Q}:=d^{1/2}{\mathfrak R}_Q={\mathfrak R}_Q \circ(-\Delta)^{1/4}-(-\Delta)^{1/4}\circ{\mathfrak R}_Q \quad.
\]
We have
\be
\label{II-f-6}
\begin{array}{l}
\ds d^{1/2}{\mathfrak R}_Q(v)(x)=\int_{\R}\int_{\R}\ {\mathcal R}_Q(x,y)\ \frac{v(y)-v(z)}{|y-z|^{3/2}}\ dy\, dz\\[5mm]
\ds\quad\quad-\int_{\R}\frac{dz}{|z-x|^{3/2}}\lf[\int_{\R}{\mathcal R}_Q(x,y)\ v(y)\ dy-\int_{\R}{\mathcal R}_Q(z,y)\ v(y)\ dy\rg]\quad.
\end{array}
\ee
Hence the corresponding Schwartz Kernel is equal to
\be
\label{II-f-7}
\begin{array}{l}
\ds{\mathcal R}_{d^{1/2}Q}(x,y):=\int_{\R} \frac{{\mathcal R}_Q(x,y)-{\mathcal R}_Q(x,z)}{|y-z|^{3/2}}\ dz+\int_{\R} \frac{{\mathcal R}_Q(y,z)-{\mathcal R}_Q(y,x)}{|x-z|^{3/2}}\ dz\quad.
\end{array}
\ee
We write ${\mathbf {y=x+h}}$ and we decompose $${\mathcal R}_{d^{1/2}Q}(x,x+h)={\mathcal R}^+_{d^{1/2}Q}(x,x+h)+{\mathcal R}^{1,+}_{d^{1/2}Q}(x,x+h)+{\mathcal R}^{1,-}_{d^{1/2}Q}(x,x+h)+{\mathcal R}^-_{d^{1/2}Q}(x,x+h)\quad,$$
where
\[
\begin{array}{l}
\ds{\mathcal R}^+_{d^{1/2}Q}(x,y):={\mathcal R}_Q(x,y)\ \int_{2\,|h|<|x-z|}  \lf[\frac{1}{|y-z|^{3/2}}-\frac{1}{|x-z|^{3/2}}\rg]\ dz\ \\[5mm] 
\ds\quad-\int_{2\,|h|<|x-z|} \frac{{\mathcal R}_Q(x,z)}{|y-z|^{3/2}}\ dz+\int_{2\,|h|<|x-z|} \frac{{\mathcal R}_Q(y,z)}{|x-z|^{3/2}}\ dz\quad,
\end{array}
\]
and
\[
\begin{array}{l}
\ds{\mathcal R}^-_{d^{1/2}Q}(x,y):=\int_{2\,|x-z|<|h|} \lf[\frac{{\mathcal R}_Q(x,y)-{\mathcal R}_Q(x,z)}{|y-z|^{3/2}}\rg]\ dz
+\int_{2\,|x-z|<|h|} \frac{{\mathcal R}_Q(y,z)-{\mathcal R}_Q(y,x)\ }{|x-z|^{3/2}}\ dz\quad,
\end{array}
\]
and
\[
\begin{array}{l}
\ds{\mathcal R}^{1,+}_{d^{1/2}Q}(x,y):=\int_{2\,|y-z|<|h|\simeq |x-z|} \lf[\frac{{\mathcal R}_Q(x,y)-{\mathcal R}_Q(x,z)}{|y-z|^{3/2}}\rg]\ dz
+\int_{2\,|y-z|<|h|\simeq |x-z|} \frac{{\mathcal R}_Q(y,z)-{\mathcal R}_Q(y,x)\ }{|x-z|^{3/2}}\ dz\quad,
\end{array}
\]
and
\[
\begin{array}{l}
\ds{\mathcal R}^{1,-}_{d^{1/2}Q}(x,y):=\int_{|y-z|\simeq |h|\simeq |x-z|} \lf[\frac{{\mathcal R}_Q(x,y)-{\mathcal R}_Q(x,z)}{|y-z|^{3/2}}\rg]\ dz
+\int_{|y-z|\simeq |h|\simeq |x-z|} \frac{{\mathcal R}_Q(y,z)-{\mathcal R}_Q(y,x)\ }{|x-z|^{3/2}}\ dz\quad.
\end{array}
\]
\par
\medskip
{\bf1.  We first bound ${\mathbf{\|{\mathcal R}^+_{d^{1/2}Q}\|_{A^{-1/2+1/p}_p}}}$}\par
\medskip
 To that aim we estimate
\be
\label{II-f-8}
\begin{array}{l}
\ds\int_{\R}h^{2-2/p}\lf|\int_{\R}|{\mathcal R}_Q(x,x+h)|^p\ \lf|\int_{2\,|h|<|x-z|}  \lf[\frac{1}{|x+h-z|^{3/2}}-\frac{1}{|x-z|^{3/2}}\rg]\ dz\rg|^pdx\rg|^{2/p}\ dh\\[5mm]
\ds\le\int_{\R}h^{2-2/p}\lf|\int_{\R}|{\mathcal R}_Q(x,x+h)|^p\ \lf|\int_{2\,|h|<|x-z|} \frac{|h|}{|x-z|^{5/2}}\ dz\rg|^p\ dx\rg|^{2/p}\ dh\\[5mm]
\ds\le\int_{\R}h^{1-2/p}\lf|\int_{\R}|{\mathcal R}_Q(x,x+h)|^p\ dx\rg|^{2/p}\ dh\quad.\\[5mm]
\end{array}
\ee
We shall now prove the following intermediate lemma
\begin{Lm}
\label{lm-interm}
Under the previous notations one has
\be
\label{II-f-4a}
\int_{\R}h^{1-2/p}\lf|\int_{\R}|{\mathcal R}_Q(x,x+h)|^p\ dx\rg|^{2/p}\ dh\le C\ \|Q\|^2_{B^{1/p}_{p,2}} \quad.
\ee
\hfill $\Box$
\end{Lm}
\noindent{\bf Proof of lemma~\ref{lm-interm}.}
\begin{equation}
\label{II-f-4b}
\begin{array}{l}
\ds\frac{1}{\pi^2}\int_{\R}h^{1-\frac{2}{p}}\lf|\int_{\R}|{\mathcal R}_Q(x,x+h)|^p\ dx\rg|^{2/p}\ dh \\[5mm]
\ds\le \frac{1}{\pi^2} \int_{\R}h^{1-\frac{2}{p}}\lf|\int_{\R}\lf|\int_{z<x+h/2}\ \frac{Q(z)-Q(x)}{(x-z)\, (z-x-h)}\ dz+\int_{z>x+h/2}\ \frac{Q(z)-Q(x+h)}{(x-z)\, (z-x-h)}\ dz\rg|^p dx\rg|^{\frac{2}{p}}\ dh 
\\[5mm]
\ds\le \frac{2}{\pi^2}\ \int_{\R}h^{1-\frac{2}{p}}\lf|\int_{\R}\lf|\int_{x-h/2<z<x+h/2}\ \frac{Q(z)-Q(x)}{(x-z)\, (z-x-h)}\ dz\rg.\rg.\\[5mm]
\ds\quad\quad\quad\quad\quad\quad\lf.\lf.+\int_{x+h-h/2<z<x+3h/2}\ \frac{Q(z)-Q(x+h)}{(x-z)\, (z-x-h)}\ dz\rg|^p dx\rg|^{\frac{2}{p}}\ dh 
 \\[5mm]
\ds+\frac{2}{\pi^2}\ \int_{\R}h^{1-\frac{2}{p}}\lf|\int_{\R}\lf|\int_{z<x-h/2}\ \frac{Q(z)-Q(x)}{(x-z)\, (z-x-h)}\ dz\rg|^p \ dx\rg|^{\frac{2}{p}}\ dh 
\\[5mm]
\ds +\frac{2}{\pi^2}\ \int_{\R}h^{1-\frac{2}{p}}\lf|\int_{\R}\lf|\int_{x+3h/2<z}\ \frac{Q(z)-Q(x+h)}{(x-z)\, (z-x-h)}\ dz\rg|^p \ dx\rg|^{\frac{2}{p}}\ dh. 
\end{array}
\end{equation}
We have on one hand, denoting $y:=x+h$
\be
\label{II-f-5b}
\begin{array}{l}
\ds\int_{\R}h^{1-\frac{2}{p}}\lf|\int_{\R}\lf|\int_{x-h/2<z<x+h/2}\ \frac{Q(z)-Q(x)}{(x-z)\, (z-y)}\ dz+\int_{y-h/2<z<y+h/2}\ \frac{Q(z)-Q(y)}{(x-z)\, (z-y)}\ dz\rg|^p \ dx\rg|^{\frac{2}{p}}\ dh \\[5mm]
\ds\le\int_{\R}h^{-1-\frac{2}{p}}\lf|\int_{\R}\lf|\int_{x-h/2<z<x+h/2}\ \frac{Q(z)-Q(x)}{(z-x)}\ dz-\int_{y-h/2<z<y+h/2}\ \frac{Q(z)-Q(x+h)}{(z-y)}\ dz\rg|^p \ dx\rg|^{\frac{2}{p}}\ dh \\[5mm]
\ds+\int_{\R}h^{-1-\frac{2}{p}}\lf|\int_{\R}\lf|\int_{x-h/2<z<x+h/2}\ \frac{|Q(z)-Q(x)|}{|h|}\ dz+\int_{y-h/2<z<y+h/2}\ \frac{|Q(z)-Q(y)|}{|h|}\ dz\rg|^p \ dx\rg|^{\frac{2}{p}}\ dh\\[5mm]
\ds\le \int_{\R}h^{-1-2/p}\lf|\int_{\R}\lf|{\mathfrak R}(Q)(x)-{\mathfrak R}(Q)(x+h)\rg|^p \ dx\rg|^{2/p}\ dh\\[5mm]
\ds+\int_{\R}h^{-1-2/p}\lf|\int_{\R} \lf|\int_{\{z<x-h/2\}\cup\{x+3h/2<z\}}\  \frac{Q(z)-Q(x)}{(z-x)}-\frac{Q(z)-Q(x+h)}{(z-x-h)}\ dz \rg|^p \ dx\rg|^{2/p}\ dh \\[5mm]
\ds+\int_{\R}h^{-1-2/p}\lf|\int_{\R}\lf|\int_{|t|<1/2}\ {|Q(x+th)-Q(x)|}\ dt+\int_{|t|<1/2}\ {|Q(y+th)-Q(y)|}\ dt\rg|^p \ dx\rg|^{2/p}\ dh\\[5mm]
\ds\le\pi^2 \|{\mathfrak R}(Q)\|^2_{B^{1/p}_{p,2}}+\int_{\R}h^{-1-2/p}\lf|\int_{-1/2}^{1/2}\|Q(x+th)-Q(x)\|_{L^p}+\|Q(x+h+th)-Q(x+h)\|_{L^p}\ dt\rg|^2\ dh\\[5mm]
\ds +\int_{\R}h^{-1-2/p}\lf|\int_{\R} \lf|\int_{\{z<x-h/2\}\cup\{y+h/2<z\}}\  h\, \frac{Q(x)-Q(z)}{(z-x)\ (z-y)}-\frac{Q(x)-Q(y)}{(z-y)}\ dz\rg|^p \ dx\rg|^{2/p}\ dh \quad.\\[5mm]
\end{array}
\ee
Using the fact that
\be
\label{II-f-6b}
\int_{\{z<x-h/2\}\cup\{x+3h/2<z\}}\ \frac{dz}{z-x-h}=\int_{h/2<u<3h/2}\frac{du}{u}=\log 3\quad,
\ee
we deduce
\be
\label{II-f-7b}
\begin{array}{l}
\ds\int_{\R}h^{1-2/p}\lf|\int_{\R}\lf|\int_{x-h/2<z<x+h/2}\ \frac{Q(z)-Q(x)}{(x-z)\, (z-y)}\ dz+\int_{y-h/2<z<y+h/2}\ \frac{Q(z)-Q(y)}{(x-z)\, (z-y)}\ dz\rg|^p \ dx\rg|^{\frac{2}{p}}\ dh \\[5mm]
\ds\le C\ \|Q\|^2_{B^{1/p}_{p,2}}+\int_{\R}h^{-1-2/p}\lf|\int_{\R} \lf|\int_{\{t<-1/2\}\cup\{3/2<t\}}\  |Q(x)-Q(x+th)|\ \frac{dt}{|t|\, |t-1|}\rg|^p \ dx\rg|^{2/p}\ dh\\[5mm]
\ds\le C\ \|Q\|^2_{B^{1/p}_{p,2}} \quad.
\end{array}
\ee
Now we treat the two last terms of the r.h.s. of (\ref{II-f-4b}). First we have
\be
\label{II-f-8b}
\begin{array}{l}
\ds\int_{\R}h^{1-2/p}\lf|\int_{\R}\lf|\int_{z<x-h/2}\ \frac{Q(z)-Q(x)}{(x-z)\, (z-x-h)}\ dz\rg|^p \ dx\rg|^{2/p}\ dh\\[5mm]
\ds\le\int_{\R}h^{-1-2/p}\lf|\int_{\R}dx\lf|\int_{t>1/2}\frac{Q(x)-Q(x+th)}{t\,  (1+t)}\ dt\rg|^p\rg|^{2/p}\ dh\\[5mm]
\ds\le\int_{\R}h^{-1-2/p}\lf|\int_{t>1/2}\frac{dt}{t\, (1+t)}\|Q(\cdot)-Q(\cdot+th)\|_p\rg|^2\ dh\\[5mm]
\ds\le C\ \int_{t>1/2}\frac{dt}{t\, (1+t)}\ \int_{\R}h^{-1-2/p}\ \|Q(\cdot)-Q(\cdot+th)\|_p^2\ dh\\[5mm]
\ds\le C\ \int_{t>1/2}\frac{t^{2/p}}{t\, (1+t)}\ dt\ \|Q\|^2_{B^{1/p}_{p,2}} = C'\ \|Q\|^2_{B^{1/p}_{p,2}}\quad.
\end{array}
\ee
We treat the  last term of the r.h.s. of (\ref{II-f-4b}) in a similar way and we establish the lemma~\ref{lm-interm}.\hfill $\Box$

\medskip


\medskip

Combining (\ref{II-f-8}) and (\ref{II-f-8b}) we have then
\be
\label{II-f-8t}
\begin{array}{l}
\ds\int_{\R}h^{2-2/p}\lf|\int_{\R}|{\mathcal R}_Q(x,x+h)|^p\ \lf|\int_{2\,|h|<|x-z|}  \lf[\frac{1}{|x+h-z|^{3/2}}-\frac{1}{|x-z|^{3/2}}\rg]\ dz\rg|^pdx\rg|^{2/p}\ dh\\[5mm]
\ds\le C\ \|Q\|^2_{B^{1/p}_{p,2}} \quad.
\end{array}
\ee
We have now
\be
\label{II-f-9}
\begin{array}{l}
\ds\int_{\R}h^{2-2/p}\lf|\int_{\R}\lf|\int_{2\,|h|<|x-z|} \frac{{\mathcal R}_Q(x,z)}{|y-z|^{3/2}}-\frac{{\mathcal R}_Q(y,z)}{|x-z|^{3/2}}\ dz\rg|^pdx\rg|^{2/p}\ dh\\[5mm]
\ds\le\ 2\ \int_{\R}h^{2-2/p}\lf|\int_{\R}\lf|\int_{2\,|h|<|x-z|}{\mathcal R}_Q(x,z)\ \lf[\frac{1}{|y-z|^{3/2}}-\frac{1}{|x-z|^{3/2}}\rg]\ dz\rg|^pdx\rg|^{2/p}\ dh\\[5mm]
\ds+\ 2\ \int_{\R}h^{2-2/p}\lf|\int_{\R}\lf|\int_{2\,|h|<|x-z|}\frac{{\mathcal R}_Q(y,z)-{\mathcal R}_Q(x,z)}{|x-z|^{3/2}}\ dz\rg|^pdx\rg|^{2/p}\ dh\quad.
\end{array}
\ee
We have first
\be
\label{II-f-10}
\begin{array}{l}
\ds\int_{\R}h^{2-2/p}\lf|\int_{\R}\lf|\int_{2\,|h|<|x-z|}{\mathcal R}_Q(x,z)\ \lf[\frac{1}{|y-z|^{3/2}}-\frac{1}{|x-z|^{3/2}}\rg]\ dz\rg|^pdx\rg|^{2/p}\ dh\\[5mm]
\ds\le\int_{\R}h^{2-2/p}\lf|\int_{\R}\lf|\int_{2\,|h|<|x-z|}|{\mathcal R}_Q(x,z)|\ \frac{|h|}{|x-z|^{5/2}}\ dz\rg|^pdx\rg|^{2/p}\ dh\\[5mm]
\ds\le\int_{\R}h^{1-2/p}\lf\|\int_2^{+\infty}\frac{dt}{t^{5/2}}\ |{\mathcal R}_Q(\cdot,\cdot+th)|\rg\|_p^2\ dh\le \int_{\R}h^{1-2/p}\lf[\int_2^{+\infty}\frac{dt}{t^{5/2}}\lf\| {\mathcal R}_Q(\cdot,\cdot+th)\rg\|_p\rg]^2\ dh\\[5mm]
\ds\le\ C\,\int_2^{+\infty}\frac{dt}{t^{5/2}}\int_{\R}h^{1-2/p}\lf\| {\mathcal R}_Q(\cdot,\cdot+th)\rg\|_p^2\ dh\\[5mm]
\ds\le C\,\int_2^{+\infty}\frac{dt}{t^{5/2+2-2/p}}\int_{\R}(th)^{1-2/p}\lf\| {\mathcal R}_Q(\cdot,\cdot+th)\rg\|_p^2\ d(th)\\[5mm]
\ds\le\ C\, \int_{\R}|h|^{1-2/p}\lf\| {\mathcal R}_Q(\cdot,\cdot+h)\rg\|_p^2\ dh\quad.
\end{array}
\ee
Combining (\ref{II-f-8})...(\ref{II-f-10})  we obtain that
\be
\label{II-f-11}
\begin{array}{l}
\ds\int_{\R}h^{2-2/p}\lf|\int_{\R}\lf|\int_{2\,|h|<|x-z|}{\mathcal R}_Q(x,z)\ \lf[\frac{1}{|y-z|^{3/2}}-\frac{1}{|x-z|^{3/2}}\rg]\ dz\rg|^pdx\rg|^{2/p}\ dh\le\  C\ \|Q\|^2_{B^{1/p}_{p,2}} \quad.
\end{array}
\ee
We have also
\be
\label{II-f-12}
\begin{array}{l}
\ds\int_{\R}h^{2-2/p}\lf|\int_{\R}\lf|\int_{2\,|h|<|x-z|}{\mathcal R}_Q(x,z)\ \frac{1}{|x-z|^{3/2}} dz\rg|^pdx\rg|^{2/p}\ dh\\[5mm]
\ds\le\int_{\R}h^{1-2/p}\lf\|\int_2^{+\infty}\frac{dt}{t^{3/2}}\ |{\mathcal R}_Q(\cdot,\cdot+th)|\rg\|_p^2\ dh\le \int_{\R}h^{1-2/p}\lf[\int_2^{+\infty}\frac{dt}{t^{3/2}}\lf\| {\mathcal R}_Q(\cdot,\cdot+th)\rg\|_p\rg]^2\ dh\\[5mm]
\ds\le\ C\,\int_2^{+\infty}\frac{dt}{t^{3/2}}\int_{\R}h^{1-2/p}\lf\| {\mathcal R}_Q(\cdot,\cdot+th)\rg\|_p^2\ dh\\[5mm]
\ds\le C\,\int_2^{+\infty}\frac{dt}{t^{3/2+2-2/p}}\int_{\R}(th)^{1-2/p}\lf\| {\mathcal R}_Q(\cdot,\cdot+th)\rg\|_p^2\ d(th)\\[5mm]
\ds\le\ C\, \int_{\R}|h|^{1-2/p}\lf\| {\mathcal R}_Q(\cdot,\cdot+h)\rg\|_p^2\ dh\quad.
\end{array}
\ee
In a similar way as in \eqref{II-f-12} one can estimate 
$$\int_{\R}h^{2-2/p}\lf|\int_{\R}\lf|\int_{2\,|h|<|x-z|}{\mathcal R}_Q(y,z)\ \frac{1}{|x-z|^{3/2}} dz\rg|^pdx\rg|^{2/p}\ dh.$$

\par
\bigskip
Combining (\ref{II-f-8})...(\ref{II-f-12}) we have proved the following lemma
\begin{Lm}
\label{lm-R^+} Under the above notations one has
\be
\label{II-f-13}
\int_{\R}h^{2-2/p}\|{\mathcal R}^+_{d^{1/2}Q}(\cdot,\cdot+h)\|_p^2\ dh\le\ C\ \|Q\|^2_{B^{1/p}_{p,2}} \quad.
\ee
\hfill $\Box$
\end{Lm}
\par
\medskip
\noindent{\bf 2. We are now bounding $\mathbf {\mathcal R}^-_{d^{1/2}Q}$. }\par
We have for $z<y$ and $x<z$
\be
\label{II-f-14}
\begin{array}{l}
\ds{\mathcal R}_Q(y,z)-{\mathcal R}_Q(y,x)=\frac{1}{\pi^2}\lf<\mbox{PV}\lf(\frac{x-z}{(\cdot-y)\, (x-\cdot)\, (z-\cdot)}\rg), Q(\cdot)\rg>\\[5mm]
\ds \quad=(x-z)\ \int_{\R}\frac{Q(\xi)}{(\xi-y)\, (x-\xi)\, (z-\xi)}\ d\xi\quad.
\end{array}
\ee
An elementary study of function gives the existence of $a\in (x,z)$ and $b\in (z,y)$ such that
\be
\label{II-f-15}
\int_{-\infty}^a\frac{1}{(\xi-y)\, (x-\xi)\, (z-\xi)}\ d\xi=0\quad\mbox{ and }\int_{a}^b\frac{1}{(\xi-y)\, (x-\xi)\, (z-\xi)}\ d\xi=0\quad,
\ee
moreover, since $2|x-z|<|z-y|$ we have in particular
\be
\label{II-f-15-a}
\min\lf\{\frac{|a-x|}{|x-z|},\frac{|a-z|}{|x-z|}\rg\}\ge \eta>0\quad\mbox{ and }\quad\min\lf\{\frac{|b-y|}{|y-z|},\frac{|b-z|}{|y-z|}\rg\}\ge \eta>0\quad,
\ee
where $\eta$ is a universal constant strictly less than $1/2$. We prove now that  $b$ is uniformly bounded from above and from below away from 0. Without loss of generality we can take $z=0$, $y=1$ and $x\in (-1/2,0)$.
We have
\[
\frac{1}{(\xi-1)\, (\xi-x)\, \xi}=\frac{1}{(\xi-1)\, x}\,\lf[\frac{1}{(\xi-x)}-\frac{1}{\xi}   \rg]=\frac{1}{x\, (x-1)}\,\lf[ \frac{1}{(\xi-x)}-\frac{1}{(\xi-1)}\rg]-\frac{1}{x}\, \lf[\frac{1}{(\xi-1)}-\frac{1}{\xi}  \rg]\quad.
\]
Hence
\[
\mbox{PV}\lf[\int_{-\infty}^t \frac{d\xi}{(\xi-1)\, (\xi-x)\, \xi}\rg]=\frac{1}{x\, (x-1)}\ \log\lf|\frac{t-x}{t-1}\rg|-\frac{1}{x}\log\lf|\frac{t-1}{t}\rg|\quad.
\]
We have in particular
\[
\frac{1}{ (1-x)}\ \log\lf|\frac{b(x)-x}{1-b(x)}\rg|=\log\lf|\frac{1-b(x)}{b(x)}\rg|\quad.
\]
This first give $b(x)<1/2$. Assume $b(x)\rightarrow 0$ as $x\rightarrow 0$ we would get
\[
0>\log(b(x)-x)\simeq \log b(x)^{-1}>0\quad,
\]
which is a contradiction.
We have then
\be
\label{II-f-16}
\begin{array}{l}
\pi^2\left(\ds{\mathcal R}_Q(y,z)-{\mathcal R}_Q(y,x)-\frac{(x-z)}{(x-y)}\ {\mathcal R}_Q(x,z)\right)=\\[5mm]
\ds=(x-z)\ \int_{-\infty}^b\frac{Q(\xi)-Q(z)}{(\xi-y)\, (x-\xi)\, (z-\xi)}-\frac{Q(\xi)-Q(z)}{(x-y)\, (x-\xi)\, (z-\xi)}\ d\xi\\[5mm]
\ds+(x-z)\ \int_{b}^{+\infty}\frac{Q(\xi)-Q(y)}{(\xi-y)\, (x-\xi)\, (z-\xi)}-\frac{Q(\xi)-Q(z)}{(x-y)\, (x-\xi)\, (z-\xi)}\ d\xi\\[5mm]
\ds=(x-z)\ \int_{-\infty}^b\frac{Q(\xi)-Q(z)}{(\xi-y)\, (x-y)\, (z-\xi)}\\[5mm]
\ds+(x-z)\ \int_{b}^{+\infty}\frac{Q(\xi)-Q(y)}{(\xi-y)\, (x-\xi)\, (z-\xi)}-(x-z)\int_{b}^{+\infty}\frac{Q(\xi)-Q(z)}{(x-y)\, (x-\xi)\, (z-\xi)}\ d\xi\quad.
\end{array}
\ee
We decompose further and we write
\be
\label{II-f-17}
\begin{array}{l}
\ds (x-z)\ \int_{-\infty}^b\frac{Q(\xi)-Q(z)}{(\xi-y)\, (x-y)\, (z-\xi)}+(x-z)\ \int_{b}^{+\infty}\frac{Q(\xi)-Q(y)}{(\xi-y)\, (x-\xi)\, (z-\xi)}\\[5mm]
\ds=(x-z)\ \int_{|\xi-z|\le \eta \,|h|}\frac{Q(\xi)-Q(z)}{(\xi-y)\, (x-y)\, (z-\xi)}+(x-z)\ \int_{|\xi-y|\le \eta\, h}\frac{Q(\xi)-Q(y)}{(\xi-y)\, (x-\xi)\, (z-\xi)}\\[5mm]
\ds+(x-z)\ \int_{z+\eta\, h}^b\frac{Q(\xi)-Q(z)}{(\xi-y)\, (x-y)\, (z-\xi)}+(x-z)\ \int_{b}^{y-\eta\, |h|}\frac{Q(\xi)-Q(y)}{(\xi-y)\, (x-\xi)\, (z-\xi)}\\[5mm]
\ds+(x-z)\ \int_{-\infty}^{z-\eta\, |h|}\frac{Q(\xi)-Q(z)}{(\xi-y)\, (x-y)\, (z-\xi)}+(x-z)\ \int_{y+\eta\,|h|}^{+\infty}\frac{Q(\xi)-Q(y)}{(\xi-y)\, (x-\xi)\, (z-\xi)}\quad.
\end{array}
\ee
We have 
\be
\label{II-f-18}
\begin{array}{l}
\ds(x-z)\ \int_{|\xi-z|\le \eta \,|h|}\frac{Q(\xi)-Q(z)}{(\xi-y)\, (x-y)\, (z-\xi)}+(x-z)\ \int_{|\xi-y|\le \eta\, |h|}\frac{Q(\xi)-Q(y)}{(\xi-y)\, (x-\xi)\, (z-\xi)}\\[5mm]
\ds+(x-z)\ \frac{{\mathfrak R}(Q)(z)-{\mathfrak R}(Q)(y)}{(x-y)\, (z-y)}\\[5mm]
\ds=\frac{(x-z)}{(x-y)}\ \int_{|\xi-z|\le \eta \,|h|}\frac{Q(\xi)-Q(z)}{(z-\xi)}\, \lf[\frac{1}{(\xi-y)}-\frac{1}{(z-y)}\rg]\ d\xi \\[5mm]
\ds+(x-z)\ \int_{|\xi-y|\le \eta \,|h|}\frac{Q(\xi)-Q(y)}{(\xi-y)}\, \lf[\frac{1}{(x-\xi)\, (z-\xi)}-\frac{1}{(x-y) \, (z-y)}\rg]\ d\xi\\[5mm]
\ds+ \frac{(x-z)}{(x-y)\, (z-y)}\int_{\{|\xi-z|>\eta\, |h|\}\cap\{|\xi-y|>\eta\, |h|\}}{[Q(\xi)-Q(z)]}\ \lf[ \frac{1}{(\xi-z)}-\frac{1}{(\xi-y)}\rg]\ d\xi\\[5mm]
\ds +\frac{(x-z)}{(x-y)\, (z-y)}\int_{|\xi-y|<\eta\, |h|}\frac{Q(\xi)-Q(z)}{(\xi-z)}\ d\xi-\frac{(x-z)}{(x-y)\, (z-y)}\int_{|\xi-z|<\eta\, |h|}\frac{Q(\xi)-Q(z)}{(\xi-y)}\ d\xi\quad.
\end{array}
\ee
\medskip

Then, we have first
\be
\label{II-f-19}
\begin{array}{l}
\ds\int_{\R_+}h^{2-2/p}\lf|\int_{\R}\lf|\int_{0<2\,(z-x)<h}\frac{1}{\sqrt{|x-z|}h}\lf|\int_{|\xi-z|\le \eta \,|h|}\frac{Q(\xi)-Q(z)}{ (z-\xi)}\ \lf[\frac{1}{(\xi-y)}-\frac{1}{(z-y)}\rg]\ d\xi\rg|\ dz\rg|^p\ dx\rg|^{2/p}\ dh\\[5mm]
\ds\le \int_{\R_+}h^{-2/p}\lf|\int_{\R}\lf|\int_{0<2\,(z-x)<h}|x-z|^{-1/2}\,\int_{|\xi-z|\le \eta \,|h|}\frac{|Q(\xi)-Q(x)|}{(z-y)^2 }\  d\xi\ dz\rg|^p\ dx\rg|^{2/p}\ dh\\[5mm]
\ds+ \int_{\R_+}h^{-2/p}\lf|\int_{\R}\lf|\int_{0<2\,(z-x)<h}|x-z|^{-1/2}\,\int_{|\xi-z|\le \eta \,|h|}\frac{|Q(z)-Q(x)|}{(z-y)^2 }\  d\xi\ dz\rg|^p\ dx\rg|^{2/p}\ dh\\[5mm]
\ds\le\int_{\R_+}h^{-4-2/p}\lf|\int_{\R}\lf|\int_{0<2\,(z-x)<h}|x-z|^{-1/2}\ dz\,\int_{|\xi-x|\le 2 \,|h|}{|Q(\xi)-Q(x)|}  d\xi\ \rg|^p\ dx\rg|^{2/p}\ dh\\[5mm]
\ds+ \int_{\R_+}h^{-2-2/p}\lf|\int_{\R}\lf|\int_{0<2\,(z-x)<h}|x-z|^{-1/2}\,{|Q(z)-Q(x)|}\  \ dz\rg|^p\ dx\rg|^{2/p}\ dh\\[5mm]
\ds\le\int_{\R_+}h^{-1-2/p}\lf|\int_{\R}\lf|\int_{0}^2\ {|Q(x+t\,h)-Q(x)|}\  dt\ \rg|^p\ dx\rg|^{2/p}\ dh\\[5mm]
\ds+ \int_{\R_+}h^{-1-2/p}\lf|\int_{\R}\lf|\int_{0}^2t^{-1/2}\,{|Q(x+t\,h)-Q(x)|}\ dt\rg|^p\ dx\rg|^{2/p}\ dh\\[5mm]
\ds\le \int_{\R_+}h^{-1-2/p}\lf|\int_{0}^2t^{-1/2}\,{\|Q(\cdot+t\,h)-Q(\cdot)\|_p}\ dt\rg|^{2}\ dh\\[5mm]
\ds\le \int_{0}^2t^{-1/2}\int_{\R_+}h^{-1-2/p}\|Q(\cdot+t\,h)-Q(\cdot)\|_p^{2}\ dh\\[5mm]
\ds\le  \int_{0}^2t^{-1/2+2/p}\int_{\R_+}h^{-1-2/p}\|Q(\cdot+h)-Q(\cdot)\|_p^{2}\ dh\le C\ \|Q\|^2_{B^{1/p}_{p,2}} \quad.
\end{array}
\ee
\par
\medskip
Similarly we have
\be
\label{II-f-20}
\begin{array}{l}
\ds\int_{\R_+}h^{2-2/p}\lf|\int_{\R}\lf|\int_{0<2\,(z-x)<h}\frac{dz}{\sqrt{|x-z|}}\,\lf|\int_{|\xi-y|\le \eta \,|h|}\frac{Q(\xi)-Q(y)}{(\xi-y)}\, \lf[\frac{1}{(x-\xi)\, (z-\xi)}\rg.\rg.\rg.\rg.\\[5mm]
\ds\lf.\lf.\lf.\lf.-\frac{1}{(x-y) \, (z-y)}\rg]\ d\xi\rg|\rg|^p\ dx\rg|^{2/p}\ dh\le C\ \|Q\|^2_{B^{1/p}_{p,2}} \quad.
\end{array}
\ee
Regarding the third term in the r.h.s. of (\ref{II-f-18}), we have, denoting $\om:=\{|\xi-z|>\eta\, |h|\}\cap\{|\xi-y|>\eta\, |h|\}$
\be
\label{II-f-21}
\begin{array}{l}
\ds\int_{\R_+}h^{2-\frac{2}{p}}\lf|\int_{\R}\lf| \int_{|x-z|<h/2}  \frac{|x-z|^{-1/2}}{|x-y|\, |z-y|}\int_{\om}{|Q(\xi)-Q(z)|}\ \lf| \frac{1}{(\xi-z)}-\frac{1}{(\xi-y)}\rg|\ d\xi\ dz\rg|^p\ dx\rg|^{\frac{2}{p}}\ dh\\[5mm]
\ds\le\int_{\R_+}h^{-2/p}\lf|\int_{\R}\lf| \int_{|x-z|<h/2} |x-z|^{-1/2}\int_{|\xi-z|>\eta\, |h|}\frac{|Q(\xi)-Q(z)|}{|\xi-z|^2}\ d\xi\ dz\rg|^p\ dx\rg|^{2/p}\ dh\\[5mm]
\ds\le\int_{\R_+}h^{-2/p}\lf|\int_{\R}\lf| \int_{|x-z|<h/2} |x-z|^{-1/2}\int_{|\xi-z|>\eta\, |h|}\frac{|Q(x)-Q(z)|}{|\xi-z|^2}\ d\xi\ dz\rg|^p\ dx\rg|^{2/p}\ dh\\[5mm]
\ds+C\,\int_{\R_+}h^{-2/p}\lf|\int_{\R}\lf| \int_{|x-z|<h/2} |x-z|^{-1/2}\int_{|\xi-x|>\eta\, |h|}\frac{|Q(\xi)-Q(x)|}{|\xi-x|^2}\ d\xi\ dz\rg|^p\ dx\rg|^{2/p}\ dh\\[5mm]
\ds\le\int_{\R_+}h^{-2-2/p}\lf|\int_{\R}\lf| \int_{|x-z|<h/2} |x-z|^{-1/2}{|Q(x)-Q(z)|} dz\rg|^p\ dx\rg|^{2/p}\ dh\\[5mm]
\ds+C\,\int_{\R_+}h^{1-2/p}\lf|\int_{\R}\lf|\int_{|\xi-x|>\eta\, |h|}\frac{|Q(\xi)-Q(x)|}{|\xi-x|^2}\ d\xi\rg|^p\ dx\rg|^{2/p}\ dh\\[5mm]
\ds\le\int_{\R_+}h^{-1-2/p}\lf|\int_{\R}\lf| \int_0^{1/2}\frac{dt}{\sqrt{t}}|Q(x)-Q(x+th)|\rg|^p\ dx\rg|^{2/p}\ dh\\[5mm]
\ds+C\,\int_{\R_+}h^{-1-2/p}\lf|\int_{\R}\lf|\int_{|t|>\eta}\frac{|Q(x+th)-Q(x)|}{t^2}\ dt\rg|^p\ dx\rg|^{2/p}\ dh\\[5mm]
\ds\le\int_{\R_+}h^{-1-\frac{2}{p}}\lf|\int_0^{1/2}\frac{dt}{\sqrt{t}}\|Q(\cdot)-Q(\cdot+th)\|_p\rg|^2\ dh\\[5mm]
\ds\quad+C\int_{\R_+}h^{-1-\frac{2}{p}}\lf|\int_\eta^{+\infty}\frac{dt}{{t}^2}\|Q(\cdot)-Q(\cdot+th)\|_p\rg|^2\ dh\\[5mm]
\ds\le C\int_0^{1/2}\frac{dt}{\sqrt{t}}\int_{\R_+}h^{-1-2/p}\|Q(\cdot)-Q(\cdot+th)\|^2_p\ dh\\[5mm]
\ds\quad+C\int_\eta^{+\infty}\frac{dt}{{t}^2}\int_{\R_+}h^{-1-2/p}\|Q(\cdot)-Q(\cdot+th)\|_p^2\ dh\le C\ \|Q\|^2_{B^{1/p}_{p,2}} \quad.
\end{array}
\ee
We estimate similarly the last term (\ref{II-f-18}). Now we estimate
\be
\label{II-f-22}
\begin{array}{l}
\ds\int_{\R_+}h^{2-2/p}\lf|\int_{\R}\lf| \int_{|x-z|<h/2}  \frac{1}{|x-z|^{1/2}\, h^2}|{\mathfrak R}(Q)(z)-{\mathfrak R}(Q)(z+h)|\ dz\rg|^p\ dx\rg|^{2/p}\ dh\\[5mm]
\ds\le\int_{\R_+}h^{-2-2/p}\lf|\int_{\R}\lf| \int_{|x-z|<h/2}  \frac{1}{|x-z|^{1/2}}|{\mathfrak R}(Q)(x)-{\mathfrak R}(Q)(z+h)|\ dz\rg|^p\ dx\rg|^{2/p}\ dh\\[5mm]
\ds+\int_{\R_+}h^{-2-2/p}\lf|\int_{\R}\lf| \int_{|x-z|<h/2}  \frac{1}{|x-z|^{1/2}}|{\mathfrak R}(Q)(z)-{\mathfrak R}(Q)(x)|\ dz\rg|^p\ dx\rg|^{2/p}\ dh\\[5mm]
\ds\le\int_{\R_+}h^{-1-2/p}\lf|\int_{\R}\lf| \int_{s<1/2}  \frac{1}{\sqrt{s}}|{\mathfrak R}(Q)(x)-{\mathfrak R}(Q)(x+(s+1)\,h)|\ ds\rg|^p\ dx\rg|^{2/p}\ dh\\[5mm]
\ds+\int_{\R_+}h^{-1-2/p}\lf|\int_{\R}\lf| \int_{s<1/2}  \frac{1}{\sqrt{s}}|{\mathfrak R}(Q)(x+s\,h)-{\mathfrak R}(Q)(x)|\ ds\rg|^p\ dx\rg|^{2/p}\ dh\\[5mm]
\le C\ \|{\mathfrak R}(Q)\|^2_{B^{1/p}_{p,2}} \le C\ \|Q\|^2_{B^{1/p}_{p,2}} \quad.
\end{array}
\ee
The 4 last terms in (\ref{II-f-17}) as well as the last integral in (\ref{II-f-16})  can be estimated using the same way as above in order to get
\be
\label{II-f-23}
\begin{array}{l}
\ds\int_{\R_+}|x-y|^{2-2/p}\lf|\int_{\R}\lf| \int_{|x-z|<|x-y|/2}  \frac{1}{|x-z|^{3/2}}\ \lf|{\mathcal R}_Q(y,z)-{\mathcal R}_Q(y,x)-\frac{(x-z)}{(x-y)}\ {\mathcal R}_Q(x,z)\rg|dz\rg|^p\ dx\rg|^{2/p}\ dy\\[5mm]
\ds\le C\ \|Q\|^2_{B^{1/p}_{p,2}} \quad.
\end{array}
\ee
Now we estimate
\be
\label{II-f-24}
\begin{array}{l}
\ds\int_{\R_+}h^{2-2/p}\lf|\int_{\R}\lf| \int_{|x-z|<h/2} \frac{{\mathcal R}_Q(x,x+h)-{\mathcal R}_Q(x,z)}{|x+h-z|^{3/2}}\ dz\rg|^p\ dx\rg|^{2/p}\ dh\quad.\\[5mm]
\end{array}
\ee
We have first
\be
\label{II-f-25}
\begin{array}{l}
\ds\int_{\R_+}h^{2-2/p}\lf|\int_{\R}\lf| \int_{|x-z|<h/2} \frac{{\mathcal R}_Q(x,z)}{|x+h-z|^{3/2}}\ dz\rg|^p\ dx\rg|^{2/p}\ dh\\[5mm]
\ds\le \int_{\R_+}h^{1-2/p}\lf|\int_{\R}\lf| \int_{0}^{1/2}{|{\mathcal R}_Q(x,x+t\,h)|}\ dt\rg|^p\ dx\rg|^{2/p}\ dh\\[5mm]
\ds\le  \int_{\R_+}h^{1-2/p}\lf|\int_{0}^{1/2}\ \|{\mathcal R}_Q(\cdot,\cdot+t\, h)\|_p\ dt\rg|^2\ dh\\[5mm]
\ds\le C_\al\ \int_{0}^{1/2}\ t^{\al}\ dt\int_{\R_+}h^{1-2/p}\|{\mathcal R}_Q(\cdot,\cdot+t\, h)\|_p^2\ dh\\[5mm]
\ds\le C_\al \ \int_{0}^{1/2}\ t^{\al-2+2/p}\ dt \int_{\R_+}h^{1-2/p}\|{\mathcal R}_Q(\cdot,\cdot+ h)\|_p^2\ dh\\[5mm]
\ds\le C_p\ \|Q\|^2_{B^{1/p}_{p,2}}\quad,
\end{array}
\ee
where we have chosen $1-2/p<\al<1$ and where we have used lemma~\ref{lm-interm}. Now we have
\be
\begin{array}{l}
\ds\int_{\R_+}h^{2-2/p}\lf|\int_{\R}\lf| \int_{|x-z|<h/2} \frac{{\mathcal R}_Q(x,x+h)}{|x+h-z|^{3/2}}\ dz\rg|^p\ dx\rg|^{2/p}\ dh\\[5mm]
\ds\le\int_{\R_+}h^{1-2/p}\|{\mathcal R}_Q(\cdot,\cdot+ h)\|_p^2\ dh\quad.
\end{array}
\ee
Finally observe that (for $p>2$)
\be
\label{II-f-24-a}
\begin{array}{l}
\ds\int_{\R_+}|x-y|^{2-2/p}\lf|\int_{\R}\lf| \int_{|x-z|>|x-y|/2}  \frac{1}{|x-z|^{3/2}}\ \lf|\frac{(x-z)}{(x-y)}\ {\mathcal R}_Q(x,z)\rg|dz\rg|^p\ dx\rg|^{2/p}\ dy\\[5mm]
\ds\le \int_{\R_+}|h|^{-2/p}\lf|\int_{\R}\lf| \int_{|x-z|>|h|/2}  \frac{\lf|{\mathcal R}_Q(x,z)\rg|}{|x-z|^{1/2}}\ dz\rg|^p\ dx\rg|^{2/p}\ dh\\[5mm]
\ds\le \int_{\R_+}|h|^{1-2/p}\lf|\int_{\R}\lf| \int_{1/2}^{+\infty}\frac{1}{\sqrt{t}}{|\mathcal R}_Q(x,x+t\, h)|\ dt\rg|^p\ dx\rg|^{2/p}\ dh\\[5mm]
\ds\le \int_{\R_+}|h|^{1-2/p}\lf|\int_{1/2}^{+\infty}\frac{1}{\sqrt{t}}\lf\|{\mathcal R}_Q(\cdot,\cdot+t\, h)\rg\|_p\ dt\rg|^{2}\ dh\\[5mm]
\ds\le\ C_\al\ \int_{1/2}^{+\infty}{{t}^\al}\int_{\R_+}|h|^{1-2/p}\lf\|{\mathcal R}_Q(\cdot,\cdot+t\, h)\rg\|^2_p\ dh\\[5mm]
\ds\le C_\al\ \int_{1/2}^{+\infty}{{t}^{\al-2+2/p}}\int_{\R_+}|h|^{1-2/p}\lf\|{\mathcal R}_Q(\cdot,\cdot+ h)\rg\|^2_p\ dh\le C_p\ \|Q\|^2_{B^{1/p}_{p,2}} \quad,
\end{array}
\ee
where we have chosen $\al>0$ such that $1-2/p<\al<1$ and where we have used lemma~\ref{lm-interm}. \par
Combining all the previous we have proved the following lemma
\begin{Lm}
\label{lm-R^-} Under the above notations one has
\be
\label{II-f-26}
\int_{\R}h^{2-2/p}\lf\|{\mathcal R}^-_{d^{1/2}Q}(\cdot,\cdot+h)-h^{-1}\ \int_{\R}{\mathcal R}_Q(\cdot,z)\ \frac{(\cdot-z)}{|\cdot-z|^{3/2}}\ dz\rg\|_p^2\ dh\le\ C\ \|Q\|^2_{B^{1/p}_{p,2}} \quad.
\ee
\hfill $\Box$
\end{Lm}
Now we bound ${\mathcal R}^{1,+}_{d^{1/2}Q}(x,y)$ and ${\mathcal R}^{1,-}_{d^{1/2}Q}(x,y)$. The only delicate term\footnote{The 3 other terms can be treated in a similar way as we did previously.} is given by
\[
\int_{2\,|y-z|< |h|\simeq |x-z|} \lf[\frac{{\mathcal R}_Q(x,y)-{\mathcal R}_Q(x,z)}{|y-z|^{3/2}}\rg]\ dz\quad.
\]
Without too much efforts one proves
\be
\label{II-f-27}
\begin{array}{l}
\ds\int_{\R}h^{2-2/p}\lf\|\int_{2\,|y-z|< |h|\simeq |x-z|} \lf[\frac{{\mathcal R}_Q(x,y)-{\mathcal R}_Q(x,z)}{|y-z|^{3/2}}\rg]\ dz\rg.\\[5mm]
\ds\quad\lf.-\int_{2\,|y-z|< |h|\simeq |x-y|} \lf[\frac{{\mathcal R}_Q(x,y)-{\mathcal R}_Q(x,z)}{|y-z|^{3/2}}\rg]\ dz\rg\|_p^2\ dh\le\ C\ \|Q\|^2_{B^{1/p}_{p,2}} \quad.
\end{array}
\ee
Now observe that the term
\[
\int_{2\,|y-z|< |h|\simeq |x-y|} \lf[\frac{{\mathcal R}_Q(x,y)-{\mathcal R}_Q(x,z)}{|y-z|^{3/2}}\rg]\ dz
\]
has already be estimated where the roles of $x$ and $y$ were exchanged and we proved while establishing lemma~\ref{lm-R^-} that $\|T_Q\|_{A_p^{-1/2+1/p}}\le  C\ \|Q\|_{B^{1/p}_{p,2}}$
where
\[
T_Q(x,y)=\int_{2\,|x-z|< |h|\simeq |x-y|} \lf[\frac{{\mathcal R}_Q(y,x)-{\mathcal R}_Q(y,z)}{|x-z|^{3/2}}\rg]\ dz+\frac{1}{(x-y)}\,\int_{\R}{\mathcal R}_Q(x,z)\ \frac{(x-z)}{|x-z|^{3/2}}\ dz\quad.
\]
Observe the following obvious fact
\be
\label{II-f-28}
\int_{\R}h^{2-2/p}\lf\|T_Q(\cdot,\cdot+h)\rg\|_p^2\ dh=\int_{\R}h^{2-2/p}\lf\|T_Q(\cdot-h,\cdot)\rg\|_p^2\ dh=\int_{\R}h^{2-2/p}\lf\|T_Q(\cdot+h,\cdot)\rg\|_p^2\ dh\quad.
\ee
The previous considerations gives that
\be
\label{II-f-29}
\lf\|\int_{2\,|y-z|< |h|\simeq |x-y|} \lf[\frac{{\mathcal R}_Q(x,y)-{\mathcal R}_Q(x,z)}{|y-z|^{3/2}}\rg]\ dz+\frac{1}{(y-x)}\,\int_{\R}{\mathcal R}_Q(y,z)\ \frac{(y-z)}{|y-z|^{3/2}}\ dz\rg\|_{A_p^{-1/2+1/p}}\le  C\ \|Q\|_{B^{1/p}_{p,2}}\quad.
\ee
Combining all the previous we obtain the following lemma
\begin{Lm}
\label{lm-d^{1/2}R}
Under the previous notations we have
\be
\label{II-f-30}
\lf\|{\mathcal R}_{d^{1/2}Q}(x,y)-\frac{1}{(x-y)}\ [F_Q(x)+F_Q(y)]\rg\|_{A_p^{-1/2+1/p}}\le  C\ \|Q\|_{B^{1/p}_{p,2}}\quad,
\ee
where 
\[
F_Q(x):=\frac{1}{\pi^2}\int_{\R}{\mathcal R}_Q(x,z)\ \frac{(x-z)}{|x-z|^{3/2}}\ dz\quad.
\]
\hfill $\Box$
\end{Lm}
We claim that $F_Q(x)=-{\mathfrak R}[(-\Delta)^{1/4}Q].$ Indeed
\begin{eqnarray}\label{est1}
F_Q(x)&:=&\frac{1}{\pi^2}\int_{\R}\int_{\R}\frac{Q(y)}{(x-y)(y-z)}\frac{(x-z)}{|x-z|^{3/2}}\ dz dy=
\frac{1}{\pi^2}\int_{\R}\int_{\R}\frac{Q(z)}{(x-z)(z-y)}\frac{(x-y)}{|x-y|^{3/2}}\ dz dy\nonumber\\[5mm]
&=&\frac{1}{\pi^2}\int_{\R}\int_{\R}\frac{Q(z)}{(x-z)(z-y)}\frac{(x-y)}{|x-y|^{3/2}}\ dz dy-\underbrace{\frac{1}{\pi^2}\int_{\R}\int_{\R}\frac{Q(x)}{(x-z)(z-y)}\frac{(x-y)}{|x-y|^{3/2}}}_{=0}\ dz dy\nonumber\\[5mm]
&+&\underbrace{\frac{1}{\pi^2}\int_{\R}\int_{\R}\frac{Q(x)-Q(y)}{(z-y)}\frac{1}{|x-y|^{3/2}}}_{=0}\ dz dy=\nonumber\\[5mm]
&=&\frac{1}{\pi^2}\int_{\R}\int_{\R}{Q(z)-Q(x)}\left[\frac{1}{(x-z)}+\frac{1}{z-y}\right]\frac{1}{|x-y|^{3/2}}\ dz dy\nonumber\\[5mm]
&+&\frac{1}{\pi^2}\int_{\R}\int_{\R}\left[\frac{(Q(x)-Q(z))}{z-y}+\frac{(Q(z)-Q(y))}{(z-y)}\right]\frac{1}{|x-y|^{3/2}}\nonumber\\[5mm]
&=&\frac{1}{\pi^2}\int_{\R}\int_{\R}{Q(z)-Q(x)}\frac{1}{(x-z)} \frac{1}{|x-y|^{3/2}}+\frac{1}{\pi^2}\int_{\R}\int_{\R}\frac{Q(z)-Q(y)}{(z-y)}\frac{1}{|x-y|^{3/2}}
\nonumber\\[5mm]
&=&\frac{1}{\pi^2}\int_{\R}\frac{1}{|x-y|^{3/2}}\left[\int_{\R}\int_{\R}{Q(z)-Q(y)}\frac{1}{(z-y)} dz-\int_{\R}{Q(z)-Q(x)}\frac{1}{(z-x)} \ dz\right]\nonumber\\[5mm]
&=&\frac{1}{\pi^2}\int_{\R}\frac{1}{|x-y|^{3/2}}\left[{\mathfrak R}Q(y)-{\mathfrak R}Q(x)\right]\, dy=-\Delta^{1/4}{\mathfrak R}(Q)(x)=-{\mathfrak R}[\Delta^{1/4}Q](x)
\end{eqnarray}

\medskip

The operator associated to ${\mathcal L}_Q(x,y):={\mathcal R}_{d^{1/2}Q}(x,y)-\frac{1}{(x-y)}\ [F_Q(x)+F_Q(y)]$ is given by
\[
{\mathfrak L}_Q:={\mathfrak R}_Q\circ(-\Delta)^{1/4}-(-\Delta)^{1/4}\circ{\mathfrak R}_Q-F_Q\circ {\mathfrak R}-{\mathfrak R}\circ F_Q\quad,
\]
where by some abuse of notation we keep denoting $F_Q$ the operator corresponding to the multiplication by $F_Q$. Observe that this operator is anti-self-dual if $Q$ is 
taking values into symmetric matrices. Hence ${\mathcal L}_Q(x,y)$ generates a multi-commutator. We compute
\be
\label{II-f-31}
\int_{\R}{\mathcal L}_Q(x,y)\ dy=\int_{\R}\int_{\R} \frac{{\mathcal R}_Q(x,y)-{\mathcal R}_Q(x,z)}{|y-z|^{3/2}}\ dz\, dy+\int_{\R}\int_{\R} \frac{{\mathcal R}_Q(y,z)-{\mathcal R}_Q(y,x)}{|x-z|^{3/2}}\ dz\, dy - \,{\mathfrak R}(F_Q)(x)\quad.
\ee
Observe that
\be
\label{II-f-32}
\begin{array}{l}
\ds\int_{\R}\int_{\R} \frac{{\mathcal R}_Q(x,y)-{\mathcal R}_Q(x,z)}{|y-z|^{3/2}}\ dz\, dy+\int_{\R}\int_{\R} \frac{{\mathcal R}_Q(y,z)-{\mathcal R}_Q(y,x)}{|x-z|^{3/2}}\ dz\, dy\\[5mm]
\ds\quad\quad=\lf[{\mathfrak R}_Q\circ(-\Delta)^{1/4}-(-\Delta)^{1/4}\circ{\mathfrak R}_Q\rg](1)=0\quad.
\end{array}
\ee
Hence we deduce the following Corollary
\begin{Co}
\label{co-II.1}
Let $Q\in \dot{H}^{1/2}({\R},\mbox{Sym}_m({\R}))$ then the operator given by
\begin{equation}\label{opL}
{\mathfrak R}\circ Q\circ{\mathfrak R}\circ(-\Delta)^{1/4}-(-\Delta)^{1/4}\circ{\mathfrak R}\circ Q\circ{\mathfrak R}-F_Q\circ {\mathfrak R}-{\mathfrak R}\circ F_Q- \,{\mathfrak R}(F_Q)
\end{equation}
is a multi-commutator in the sense that it sends $L^2({\R})$ into $B^{-1/2+1/p}_{2p/(p+2),2}({\R})\hookrightarrow H^{-1/2}({\R}).$
\hfill $\Box$
\end{Co}
 
\begin{Rm}
We knew that 
\be\label{comm1}
\begin{array}{l}
\ds{\mathfrak R}\circ Q\circ{\mathfrak R}\circ(-\Delta)^{1/4}-(-\Delta)^{1/4}\circ{\mathfrak R}\circ Q\circ{\mathfrak R}-{\mathfrak R}\circ (-\Delta)^{1/4}Q\circ {\mathfrak R}\\[5mm]
\ds={\mathfrak R}\circ\lf[ Q\circ (-\Delta)^{1/4}-(-\Delta)^{1/4}\circ Q-(-\Delta)^{1/4}(Q)\rg]\circ{\mathfrak R}
\end{array}
\ee
and 
\be\label{comm2}
{\mathfrak R}\circ\lf[ (-\Delta)^{1/4}Q\circ {\mathfrak R}+{\mathfrak R}((-\Delta)^{1/4}Q)\rg]
\ee
are enjoying compensation property (the first one (\ref{comm1}) is the adjoint action of the Riesz transform on a 3-commutator, the second one (\ref{comm2}) is the composition between Riesz and a Coifman-Rochberg-Weiss commutator)  .
If we sum \eqref{comm1} and \eqref{comm2} we deduce that
 
\[
{\mathfrak R}\circ Q\circ{\mathfrak R}\circ(-\Delta)^{1/4}-(-\Delta)^{1/4}\circ{\mathfrak R}\circ Q\circ{\mathfrak R}-{\mathfrak R}\circ F_Q
\]
has compensation properties.
Since $F_Q=-{\mathfrak R}((-\Delta)^{1/4}Q)\in L^2$ we also have 
\[
-F_Q\circ {\mathfrak R}- \,{\mathfrak R}(F_Q)
\]
is a Coifman-Rochberg-Weiss commutator. Hence it was known that ${\mathcal T}_{{\mathcal L}_Q}$ enjoyed compensation properties. The novelty in corollary~\ref{co-II.1} is the estimate of the kernel  ${\mathcal L}_Q$  in  $A^{-1/2+1/p}_p$ which makes it a multi-commutator enjoying stability by the adjoint actions of elements in $L^\infty({\R},M_m({\R}))$ for instance.
\hfill $\Box$
\end{Rm}

\end{document}